\documentclass[10pt, a4paper]{article}
\usepackage{amscd}
\usepackage{amsfonts}
\begin{document}

\newtheorem{theorem}{Theorem}
\newtheorem{corollary}[theorem]{Corollary}
\newtheorem{definition}[theorem]{Definition}
\newtheorem{lemma}[theorem]{Lemma}
\newtheorem{proposition}[theorem]{Proposition}
\newtheorem{remark}[theorem]{Remark}
\newtheorem{example}[theorem]{Example}
\newtheorem{notation}[theorem]{Notation}
\def\Qed{\hfill\raisebox{.6ex}{\framebox[2.5mm]{}}\\[.15in]}

\title{Involutions on surfaces with $p_g=q=1$}

\author{Carlos Rito}

\date{}
\pagestyle{myheadings}
\maketitle
\setcounter{page}{1}

\begin{abstract}

In this paper some numerical restrictions for surfaces with an involution are obtained.
These formulas are used to study surfaces of general type $S$ with $p_g=q=1$
having an involution $i$ such that $S/i$ is a non-ruled surface
and such that the bicanonical map of $S$ is not composed with $i.$
A complete list of possibilities is given and several new examples are constructed,
as bidouble covers of surfaces. In particular the first example of a minimal surface of general type with
$p_g=q=1$ and $K^2=7$ having birational bicanonical map is obtained.

\noindent 2000 Mathematics Classification: 14J29.
\end{abstract}

\tableofcontents

\section{Introduction}

Several authors have studied surfaces of general type with $p_g=q=1$
(\cite{Ca1}, \cite{Ca2}, \cite{CC1}, \cite{CC2}, \cite{CP},
\cite{Po1}, \cite{Po2}, \cite{Po3}, \cite{Pi}),
but these surfaces are still not completely understood.

In \cite{Ri2}, the author gives several new examples of double planes of general type with $p_g=q=1$
having bicanonical map $\phi_2$ composed with the corresponding involution.
The case $S/i$ non-ruled and $\phi_2$ composed with $i$ is considered in \cite{Ri1}.

In this paper we study the case $\phi_2$ not composed with $i.$ More precisely,
we consider surfaces of general type $S$ with $p_g=q=1$ having an involution $i$
such that $S/i$ is a non-ruled surface and $\phi_2$ is not composed with $i.$
A list of possibilities is given and new examples are obtained for each value of
the birational invariants of $S/i$ (only the existence of the case
${\rm Kod}(S/i)=2,$ $\chi(S/i)=1$ and $q(S/i)=0$ remains an open problem).

The paper is organized as follows. In Section \ref{Involution} we obtain
formulas which, for a surface $S$ with an involution $i,$ relate the invariants
of $S$ and $S/i$ with the branch locus of the cover $S\rightarrow S/i,$
its singularities and the number of nodes of $S/i.$
Section \ref{q1} contains a description of the action of the
involution $i$ on the Albanese fibration of $S.$
In Section \ref{Chapternotcomposed} we apply the numerical formulas of
Section \ref{Involution} to the case $p_g=q=1,$ obtaining the classification results,
Theorems \ref{prop0}, \ref{prop1} and \ref{prp2}. Results of Miyaoka and Sakai
on the maximal number of disjoint smooth rational or elliptic curves on a surface
are also used here. Finally Section \ref{bdce} contains the construction of examples,
as bidouble covers of surfaces:\\
the surfaces constructed in Sections \ref{expl4}, \ref{expl5} and \ref{expl3}
are Du Val double planes ({\it cf.} \cite{Ri2}) which have other interesting involutions;\\
Section \ref{expl6} contains the construction of a surface with $K^2=4,$ Albanese fibration of genus $g=2$
and ${\rm deg}(\phi_2)=2$ (thus it is not the example in \cite{Ca2}, for which
$\phi_2$ is composed with the three involutions associated to the bidouble cover);\\
in Section \ref{expl8} a new surface with $K^2=8$ is obtained (it is not
a standard isotrivial fibration);\\
Sections \ref{expl9} and \ref{expl10} contain the construction of new surfaces with
$K^2=7,6$ and ${\rm deg}(\phi_2)=1$ (it is the first example with $p_g=q=1$ and $K^2=7$
having birational bicanonical map);\\
bidouble covers of irregular ruled surfaces give interesting examples in
Sections \ref{expl7} and \ref{expl11}.

Some branch curves for these bidouble cover examples are computed in Appendix \ref{CuM},
using the {\em Computational Algebra System Magma} (\cite{BCP}).

\bigskip
\noindent{\bf Notation and conventions}

We work over the complex numbers;
all varieties are assumed to be projective algebraic.
For a projective smooth surface $S,$ the {\em canonical class} is denoted by
$K,$ the {\em geometric genus} by $p_g:=h^0(S,\mathcal O_S(K)),$ the {\em irregularity}
by $q:=h^1(S,\mathcal O_S(K))$ and the {\em Euler characteristic} by
$\chi=\chi(\mathcal O_S)=1+p_g-q.$

An {\em $(-n)$-curve} $C$ on a surface is a curve
isomorphic to $\mathbb P^1$ such that $C^2=-n$. We say that a curve
singularity is {\em negligible} if it is either a double point or a
triple point which resolves to at most a double point after one
blow-up.
An $(m_1,m_2,\ldots)$-point, or point of order $(m_1,m_2,\ldots),$
is a point of multiplicity $m_1,$ which resolves to a point
of multiplicity $m_2$ after one blow-up, etc.

An {\em involution} of a surface $S$ is an
automorphism of $S$ of order 2. We say that a map is {\em composed with
an involution} $i$ of $S$ if it factors through the double cover $S\rightarrow
S/i.$

The rest of the notation is standard in Algebraic Geometry.

\bigskip
\noindent{\bf Acknowledgements}

The author is a collaborator of the Center for Mathematical
Analysis, Geometry and Dynamical Systems of Instituto Superior T\'
ecnico, Universidade T\' ecnica de Lisboa, and
is a member of the Mathematics Department of the
Universidade de Tr\'as-os-Montes e Alto Douro.
This research was partially supported by FCT (Portugal) through
Project POCTI/MAT/44068/2002.

\section{Results on involutions}\label{ChptrInv}

\subsection{General facts}

Let $S$ be a smooth minimal surface of general type with an
involution $i.$ Since $S$ is minimal of general type, this
involution is biregular. The fixed locus of $i$ is the union of a
smooth curve $R''$ (possibly empty) and of $t\geq 0$ isolated points
$P_1,\ldots,P_t.$ Let $S/i$ be the quotient of $S$ by $i$ and
$p:S\rightarrow S/i$ be the projection onto the quotient. The
surface $S/i$ has nodes at the points $Q_i:=p(P_i),$ $i=1,\ldots,t,$
and is smooth elsewhere. If $R''\not=\emptyset,$ the image via $p$
of $R''$ is a smooth curve $B''$ not containing the singular points
$Q_i,$ $i=1,\ldots,t.$ Let now $h:V\rightarrow S$ be the blow-up of
$S$ at $P_1,\ldots,P_t$ and set $R'=h^*(R'').$ The involution $i$
induces a biregular involution $\widetilde{i}$ on $V$ whose fixed
locus is $R:=R'+\sum_1^t h^{-1}(P_i).$ The quotient
$W:=V/\widetilde{i}$ is smooth and one has a commutative diagram:
$$
\begin{CD}\ V@>h>>S\\ @V\pi VV  @VV p V\\ W@>g >> S/i
\end{CD}
$$
where $\pi:V\rightarrow W$ is the projection onto the quotient and
$g:W\rightarrow S/i$ is the minimal desingularization map. Notice
that $$A_i:=g^{-1}(Q_i),\ \ i=1,\ldots,t,$$ are $(-2)$-curves and
$\pi^*(A_i)=2\cdot h^{-1}(P_i).$ 

Set $B':=g^*(B'').$ Since $\pi$ is
a double cover with branch locus $B'+\sum_1^t A_i,$ it is determined
by a line bundle $L$ on $W$ such that $$2L\equiv B:=B'+\sum_1^t A_i.$$
It is well known that ({\it cf.} \cite[Chapter V, Section 22]{BPV}):\\
\begin{equation}\label{relate0}
\begin{array}{c}
  p_g(S)=p_g(V)=p_g(W)+h^0(W,\mathcal{O}_W(K_W+L)), \\\\
  q(S)=q(V)=q(W)+h^1(W,\mathcal{O}_W(K_W+L))
\end{array}
\end{equation}
and
\begin{equation}\label{relate}
\begin{array}{c}
  K_S^2-t=K_V^2=2(K_W+L)^2, \\\\
  \chi(\mathcal{O}_S)=\chi(\mathcal{O}_V)=2\chi(\mathcal{O}_W)+\frac{1}{2}L(K_W+L).
\end{array}
\end{equation}
\\
Denote by $\phi_2$ the bicanonical map of $S$ (given by $|2K|$).
From the papers \cite{CM} and \cite{CCM},
$$
 \phi_2 {\rm\ is\ composed\ with\ } i {\rm\ if\ and\ only\ if\ } h^0(W,\mathcal{O}_W(2K_W+L))=0.
$$

\subsection{Numerical restrictions}\label{Involution}

Let $P$ be a minimal model of the resolution $W$ of $S/i$ and
$\rho:W\rightarrow P$ be the corresponding projection. Denote by
$\overline{B}$ the projection $\rho(B)$ and by $\delta$ the
"projection" of $L.$
\begin{remark}\label{CanRes}
If $\overline B$ is singular, there are exceptional divisors $E_i$ and numbers $r_i\in 2\mathbb N$ such that
$$
\begin{array}{l}
E_i^2=-1,\\
K_W\equiv\rho^*(K_P)+\sum E_i,\\
2L\equiv B=\rho^*(\overline{B})-\sum r_iE_i\equiv\rho^*(2\delta)-\sum r_iE_i.
\end{array}
$$
\end{remark}
\begin{proposition}\label{equation}
With the previous notation, if $S$ is a surface of general type then:
\begin{description}
  \item[a)] $\chi(\mathcal{O}_P)-\chi(\mathcal{O}_S)=
  K_P(K_P+\delta)+\frac{1}{2}\sum(r_i-2)-h^0(W,\mathcal{O}_W(2K_W+L));$
  \item[b)] $\delta^2=-2\chi(\mathcal{O}_P)-2K_P^2-3K_P\delta+$\\
$+\frac{1}{4}\sum(r_i-2)(r_i-4)+2h^0(W,\mathcal{O}_W(2K_W+L)).$
\end{description}
\end{proposition}

\begin{proposition}\label{numbernodes}
Let $t$ be the number of nodes of $S/i.$ One has:
 \begin{description}
  \item[a)] $t=K_S^2+6\chi(\mathcal{O}_W)-2\chi(\mathcal{O}_S)-2h^0(W,\mathcal{O}_W(2K_W+L));$
  \item[b)] $t=K_SR''+8\chi(\mathcal{O}_W)-4\chi(\mathcal{O}_S)\geq 8\chi(\mathcal O_W)-4\chi(\mathcal O_S);$
  \item[c)] $K_S^2\geq 2\chi(\mathcal O_W)-2\chi(\mathcal O_S)+2h^0(W,\mathcal{O}_W(2K_W+L)).$
 \end{description}
\end{proposition}

\begin{proposition}\label{m}
With the above notation:
\begin{description}
  \item[a)] $h^0(W,\mathcal{O}_W(2K_W+L))\leq
  \frac{1}{3}K_W^2-\chi(\mathcal{O}_W)+\frac{11}{3}\chi(\mathcal{O}_S)+\frac{1}{27}K_S^2;$
  \item[b)] $h^0(W,\mathcal{O}_W(2K_W+L))\leq
  \frac{1}{2}K_W^2+5\chi(\mathcal{O}_S)+2q(S)-3\chi(\mathcal{O}_W)-2q(W).$
\end{description}
\end{proposition}
{\bf Proof of Proposition \ref{equation}:}\ ({\it cf.} \cite{CM})\\\\
a)
  From the Kawamata-Viehweg's vanishing theorem (see {\it e.g.}
  \cite[Corollary 5.12, c)]{EV}), one has
  $$h^i(W,\mathcal{O}_W(2K_W+L))=0,\ i=1,2.$$
  The Riemann-Roch theorem implies
  $$\chi(\mathcal{O}_W(2K_W+L))=\chi(\mathcal{O}_W)+\frac{1}{2}L(K_W+L)+K_W(K_W+L),$$
  thus, using (\ref{relate}),
  \begin{equation}\label{Mrg2}
  h^0(W,\mathcal{O}_W(2K_W+L))=\chi(\mathcal O_S)-\chi(\mathcal O_W)+K_W(K_W+L).
  \end{equation}
  With the notation of Remark \ref{CanRes}, we can write
  $$\chi(\mathcal{O}_P)-\chi(\mathcal{O}_S)=\frac{1}{2}K_W(2K_W+2L)-h^0(W,\mathcal{O}_W(2K_W+L))=$$
  $$=\frac{1}{2}\left(\rho^*(K_P)+\sum E_i\right)
  \left(2\rho^*(K_P+\delta)+\sum(2-r_i)E_i\right)-h^0(W,\mathcal{O}_W(2K_W+L))=$$
  $$=K_P(K_P+\delta)+\frac{1}{2}\sum(r_i-2)-h^0(W,\mathcal{O}_W(2K_W+L)).$$
\\\\
b)
From the proof of a),
$$h^0(W,\mathcal{O}_W(2K_W+L))=\chi(\mathcal{O}_W(2K_W+L))=\chi(\mathcal{O}_W)+\frac{1}{2}(2K_W+L)(K_W+L).$$
Using Remark \ref{CanRes} this means
$$h^0(W,\mathcal{O}_W(2K_W+L))=\chi(\mathcal{O}_P)+$$
$$+\frac{1}{2}\left(\rho^*(2K_P+\delta)+\frac{1}{2}\sum (4-r_i)E_i\right)\left(\rho^*(K_P+\delta)+\frac{1}{2}\sum(2-r_i)E_i\right)=$$
$$=\chi(\mathcal{O}_P)+K_P^2+\frac{3}{2}K_P\delta+\frac{1}{2}\delta^2-\frac{1}{8}\sum(r_i-2)(r_i-4).$$
\Qed
\\
{\bf Proof of Proposition \ref{numbernodes}:}\\\\
a) From formulas (\ref{relate}) and (\ref{Mrg2}),
$$t=K_S^2-2K_W(K_W+L)-2L(K_W+L)=$$
$$=K_S^2+2\chi(\mathcal{O}_S)-2\chi(\mathcal{O}_W)-2h^0(W,\mathcal{O}_W(2K_W+L))
-4\chi(\mathcal{O}_S)+8\chi(\mathcal{O}_W).$$
\\\\
b) (This is also a consequence of the holomorphic fixed point formula.)\\
From (\ref{relate}),
$$4\chi(\mathcal{O}_S)-8\chi(\mathcal{O}_W)=2L(K_W+L)=\left(B'+\sum_1^tA_i\right)(K_W+L)=$$
$$=B'(K_W+L)-t=\frac{1}{2}\pi^*(B')\pi^*(K_W+L)-t=R''K_S-t.$$
Since $S$ is of general type, $K_SR''\geq 0,$ thus $$t\geq 8\chi(\mathcal{O}_W)-4\chi(\mathcal{O}_S).$$
\\\\
c) This is immediate from a) and b).\Qed
\\
{\bf Proof of Proposition \ref{m}:}\\\\
a) This inequality is given by the following three claims.\\\\
    {\bf Claim 1:} \\
    {\em $1-p_a(B')=3\chi(\mathcal{O}_W)-3\chi(\mathcal{O}_S)-K_S^2-K_W^2+3h^0(W,\mathcal{O}_W(2K_W+L)).$}\\\\
    {\em Proof}\ : Formulas (\ref{relate}) and (\ref{Mrg2}) give
    $$L^2-K_W^2=$$
    $$=[2\chi(\mathcal O_S)-4\chi(\mathcal O_W)-LK_W]-
    [h^0(W,\mathcal{O}_W(2K_W+L))-\chi(\mathcal O_S)+\chi(\mathcal O_W)-K_WL],$$
    thus
    \begin{equation}\label{x}
    L^2=K_W^2+3\chi(\mathcal{O}_S)-5\chi(\mathcal{O}_W)-h^0(W,\mathcal{O}_W(2K_W+L)).
    \end{equation}
    Now we perform a straightforward calculation using the adjunction formula,
    (\ref{relate}), Proposition \ref{numbernodes}, a) and (\ref{x}):
    $$2p_a(B')-2=$$
    $$=K_WB'+B'^2=K_W2L+(2L)^2+2t=2L(K_W+L)+2t+2L^2=$$
    $$2[2\chi(\mathcal{O}_S)-4\chi(\mathcal{O}_W)]+$$
    $$+2[K_S^2+6\chi(\mathcal{O}_W)-2\chi(\mathcal{O}_S)-2h^0(W,\mathcal{O}_W(2K_W+L))]+$$
    $$+2[K_W^2+3\chi(\mathcal{O}_S)-5\chi(\mathcal{O}_W)-h^0(W,\mathcal{O}_W(2K_W+L))]=$$
    $$=2K_S^2+2K_W^2+6\chi(\mathcal{O}_S)-6\chi(\mathcal{O}_W)-6h^0(W,\mathcal{O}_W(2K_W+L)).\ \diamondsuit$$\\\\
Denote by $\tau$ the number of rational curves of $B'.$\\
    {\bf Claim 2:} {\em $$1-p_a(B')\leq\tau.$$}\\\\
    {\em Proof}\ :
    Write $$B'=\sum_1^{\tau} B_i'+\sum_{\tau+1}^hB_i'$$
    as a decomposition of $B'$ in (smooth) connected components such that
    $B_i',$ $i\leq\tau,$ are the rational ones.
    The adjunction formula gives
    $$2p_a(B')-2=\sum_1^h\left(K_WB_i'+B_i'^2\right)=
    \sum_1^{\tau}(2g(B_i')-2)+\sum_{\tau+1}^h(2g(B_i')-2)\geq -2\tau.\ \diamondsuit$$\\
    {\bf Claim 3:} {\em $$\tau\leq 8\left(\chi(\mathcal{O}_S)-\frac{1}{9}K_S^2\right).$$}
     \\\\
    {\em Proof}\ : Since $B'$ does not contain $(-2)$-curves and it is contained in the branch locus of
    the cover $\pi:V\rightarrow W,$ then each rational curve in $B'$ corresponds to a rational curve in $S.$
    Now the result follows from Proposition \ref{Miy} below.\ $\diamondsuit$\\\\
Therefore $1-p_a(B')\leq
8\left(\chi(\mathcal{O}_S)-\frac{1}{9}K_S^2\right)$ and using Claim
1 we obtain
the desired inequality.\\\\
b) Proposition \ref{numbernodes}, a) says that
$$K_V^2=K_S^2-t=2\chi(\mathcal{O}_S)-6\chi(\mathcal{O}_W)+2h^0(W,\mathcal{O}_W(2K_W+L)).$$
The second Betti number $b_2$ of a surface $X$ satisfies
$$b_2(X)=12\chi(\mathcal{O}_X)-K_X^2+4q(X)-2.$$
Therefore
$$b_2(V)=10\chi(\mathcal{O}_V)+6\chi(\mathcal{O}_W)+4q(V)-2-2h^0(W,\mathcal{O}_W(2K_W+L)).$$
Since $b_2(V)\geq b_2(W),$ one has the result.\Qed

\begin{proposition}\label{Miy}
{\em (\cite[Proposition 2.1.1]{Mi})}
Let $X$ be a minimal surface of non-negative Kodaira dimension.
Then the number of disjoint smooth rational curves in $X$ is bounded by
$$8\left(\chi(\mathcal O_X)-\frac{1}{9}K_X^2\right).$$
\end{proposition}

\subsection{Surfaces with an involution and $q=1$}\label{q1}
Let $S$ be a surface of general type with $q=1.$ Then the Albanese
variety of $S$ is an elliptic curve $E$ and the Albanese map is a
connected fibration (see {\it e.g.} \cite{Be} or \cite{BPV}).

Suppose that $S$ has an involution $i.$ Then $i$ preserves the
Albanese fibration (because $q(S)=1$) and so we have a commutative
diagram
\begin{equation}\label{CmtDgr}
\begin{CD}\ V@>h>>S @>>>E\\ @V\pi VV  @VV p V @VV V\\ W@>>> S/i @>>>\Delta
\end{CD}
\end{equation}
where $\Delta$ is a curve of genus $\leq 1.$ Denote by
$$f_A:W\rightarrow\Delta$$ the fibration induced by the Albanese
fibration of $S.$

Recall that $$\rho:W\rightarrow P$$ is the projection of $W$ onto
its minimal model $P$ and $$\overline{B}:=\rho(B),$$ where
$B:=B'+\sum_1^t A_i\subset W$ is the branch locus of $\pi.$
Let $$\overline{B'}:=\rho(B')\ \ \ {\rm and}\ \ \ \overline{A_i}=\rho(A_i).$$

When $\overline{B}$ has only negligible singularities the map $\rho$
contracts only exceptional curves contained in fibres of $f_A.$ In
fact otherwise there exists a $(-1)$-curve $J\subset W$ such that
$JB=2$ and $\pi^*(J)$ is transverse to the
fibres of the (genus 1 base) Albanese fibration of $S.$ This is
impossible because $\pi^*(J)$ is a rational curve.
Moreover $\rho$ contracts no curve meeting $\sum A_i,$
thus the singularities of $\overline{B}$ are exactly the
singularities of $\overline{B'},$ {\it i.e.} $\overline{B'}\bigcap\sum
\overline{A_i}=\emptyset.$ We denote the image of $f_A$ on $P$
by $\overline{f_A}.$

If $\Delta\cong\mathbb P^1$ then the double cover $E\rightarrow\Delta$ is
ramified over 4 points $p_j$ of $\Delta,$ thus the branch locus
$B'+\sum_1^t A_i$ is contained in 4 fibres $$F_A^j:=f_A^*(p_j),\
j=1,...,4,$$ of the fibration $f_A$.
Hence by Zariski's Lemma (see {\it e.g.} \cite{BPV}) the irreducible
components $B_i'$ of $B'$ satisfy $B_i'^2\leq 0.$ If $\overline{B}$
has only negligible singularities then also $\overline{B'}^2\leq 0.$
Since  $\pi^*(F_A^j)$ has even multiplicity, each component of
$F_A^j$ which is not a component of the branch locus $B'+\sum_1^t
A_i$ must be of even multiplicity.

\section{Classification results}\label{Chapternotcomposed}

From now on $S$ is a smooth minimal surface of general type with $p_g=q=1$ having an
involution $i$ such that the bicanonical map $\phi_2$ of $S$ is not composed with $i.$
Notice that then $2\leq K_S^2\leq 9,$ by the Debarre's inequality for an irregular surface
$\left( K_S^2\geq 2p_g\right) $ and by the Miyaoka-Yau inequality
$\left( K_S^2\leq 9\chi(\mathcal O_S)\right) $.

Recall from Section \ref{Involution} that $$h^0(W,\mathcal{O}_W(2K_W+L))\ne 0,$$ where
$W$ is the minimal resolution of $S/i$ and $L\equiv\frac{1}{2}B$ is
the line bundle which determines the double cover $V\rightarrow W.$

Let $P$ be a minimal model of $W$ and $\delta,$ $\overline{B}\equiv 2\delta$
and the numbers $r_i$ be as defined in Section \ref{Involution}.
Recall that $t$ denotes the number of nodes of $S/i.$
Notice that $p_g(P)\leq p_g(S)=1$ and $q(P)\leq q(S)=1.$

In the next sections the following result is useful:

\begin{proposition}\label{Sa}
{\em (\cite{Sa})}
Let $S$ be a minimal smooth surface of general type and $C\subset S$ be a disjoint union
of smooth elliptic curves. Then $$-C^2\leq 36\chi(\mathcal O_S)-4K_S^2.$$
\end{proposition}
{\bf Proof:} This follows from the inequality of \cite[Corollary 7.8]{Sa},
using\\ $KC+C^2=2p_a(C)-2=0.$ \Qed

\subsection{The case ${\rm Kod}(S/i)=0$}\label{secnotcomp}

Here we give a list of possibilities for the case ${\rm Kod}(S/i)=0.$

\begin{theorem}\label{prop0}

Let $S$ and $P$ be as above.
If ${\rm Kod}(P)=0,$ only the following cases can occur:
\begin{description}
  \item[a)] $P$ is an Enriques surface and
  \begin{description}
    \item[$\cdot$] $\{r_i\ne2\}=\{4\},$ $\overline B^2=0,$ $t-2=K_S^2\in\{2,\ldots,7\},$ or
    \item[$\cdot$] $\{r_i\ne2\}=\{4,4\},$ $\overline B^2=8,$ $t=K_S^2\in\{4,\ldots,8\},$ or
    \item[$\cdot$] $\{r_i\ne2\}=\{6\},$ $\overline B^2=16,$ $t=K_S^2\in\{4,\ldots,8\};$
  \end{description}

  \item[b)] $P$ is a bielliptic surface and
  \begin{description}
    \item[$\cdot$] $\{r_i\ne2\}=\emptyset,$ $\overline B^2=8,$ $t=0,$ $K_S^2=4,$ or
    \item[$\cdot$] $\{r_i\ne2\}=\{4\},$ $\overline B^2=16,$ $t+6=K_S^2=6$ or $7,$ or
    \item[$\cdot$] $\{r_i\ne2\}=\{4,4\},$ $\overline B^2=24,$ $t=0,$ $K_S^2=8,$ or
    \item[$\cdot$] $\{r_i\ne2\}=\{6\},$ $\overline B^2=32,$ $t=0,$ $K_S^2=8.$
  \end{description}
\end{description}
Furthermore, there are examples for
\begin{description}
  \item[$\cdot$] {\rm a)} with $K_S^2=8;$
  \item[$\cdot$] {\rm b)} with $K_S^2=4,$ $6,$ $7$ or $8.$
\end{description}
\end{theorem}
{\bf Proof :}
It is easy to see that $P$ cannot be a $K3$ surface: in this case
we get from Proposition \ref{m}, b) that
$$K_W^2\geq 2h^0(W,\mathcal{O}_W(2K_W+L))-2,$$
which implies $h^0(W,\mathcal{O}_W(2K_W+L))=1$ and $K_W^2=0.$
This contradicts the fact $\sum(r_i-2)=4\ne 0,$
given by Proposition \ref{equation}, a).

So, from the classification of surfaces (see {\it e.g.} \cite{Be} or \cite{BPV}),
$p_g(P)=q(P)=0$ or $p_g(P)=0,$ $q(P)=1$ (notice that $p_g(P), q(P)\leq 1$),
{\it i.e.} $P$ is an Enriques surface or a bielliptic surface.
\begin{description}
  \item[a)] Suppose $P$ is an Enriques surface: Proposition \ref{m}, a) implies
  that\\ 
  $h^0(W,\mathcal{O}_W(2K_W+L))\leq 3,$ with equality holding only if $K_W^2=0.$
  In this case the branch locus $\overline B$ is smooth, {\it i.e.} $\sum(r_i-2)=0,$
  which contradicts Proposition \ref{equation}, a).
  Therefore {$h^0(W,\mathcal{O}_W(2K_W+L))=1$} or $2.$

  Now the only possibilities allowed by Propositions \ref{equation}
  and \ref{numbernodes}, a), b) are:
  \begin{description}
    \item[1)] $\sum(r_i-2)=2,$ $\overline{B}^2=0,$ $t=K_S^2+2\geq 4;$
    \item[2)] $\sum(r_i-2)=4,$ $\overline{B}^2=8$ or 16, $t=K_S^2\geq 4.$
  \end{description}
  Moreover, if a nodal curve $A_i\subset B$ is not contracted to a point,
  then it is mapped onto a nodal curve of the Enriques surface $P.$
  Indeed, from the adjunction formula, $K_WA_i=0,$
  which means that $A_i$ does not intersect any $(-1)$-curve of $W.$

  An Enriques surface has at most 8 disjoint $(-2)$-curves. In case 1),
  the non-negligible singularities of $\overline B$ are a $4$-uple or $(3,3)$-point,
  hence $t\leq 9.$
  In case 2), $t=9$ only if $\overline B$ has a $(3,3)$-point, which implies
  that $S$ has an elliptic curve with negative self-intersection. Since in this case
  $K_S^2=9,$ this is impossible from Proposition \ref{Sa},
  therefore $t\leq 8.$
  \item[b)] Suppose $P$ is a bielliptic surface: from Proposition \ref{m}, a), one has\\
  $h^0(W,\mathcal{O}_W(2K_W+L))\leq 4,$ with equality holding only if $K_W^2=0.$
  In this case we get from Proposition \ref{equation}, a)
  that $$\sum(r_i-2)=2h^0(W,\mathcal{O}_W(2K_W+L))-2=6\ne 0,$$ which
  contradicts $K_W^2=0.$ Hence $h^0(W,\mathcal{O}_W(2K_W+L))\leq 3.$

  As in a), if a $(-2)$-curve $A_i\subset B$ is not contracted to a point, then it
  is mapped onto a $(-2)$-curve of $P.$ But a bielliptic surface has no $(-2)$-curves
  (from Proposition \ref{Miy}), thus the nodal curves of $B$ are contracted to singularities of $\overline B.$

  Using Propositions \ref{equation} and \ref{numbernodes}, a) one obtains the following possibilities:
    \begin{description}
    \item[1)] $\sum(r_i-2)=0,$ $\overline{B}^2=8,$ $K_S^2=t+4;$
    \item[2)] $\sum(r_i-2)=2,$ $\overline{B}^2=16,$ $K_S^2=t+6;$
    \item[3)] $\sum(r_i-2)=4,$ $\overline{B}^2=24,$ $K_S^2=t+8;$
    \item[4)] $\sum(r_i-2)=4,$ $\overline{B}^2=32,$ $K_S^2=t+8.$
  \end{description}
In case 1), $t=0,$ because $\overline B$ has only negligible singularities.
In case 2), $\overline B$ can have a $(3,3)$-point, thus
$t=0$ or $1.$ In case 3), $t=1$ only if $\overline B$ has a $(3,3)$-point, but then
$K_S^2=9$ and $S$ has an elliptic curve, which is impossible from Proposition \ref{Sa}.
Finally, in case 4), the only non-negligible singularity of  $\overline B$
is a point of multiplicity $6$ (from Proposition \ref{equation}, b)), thus $t=0.$
\end{description}
The examples are constructed in Sections \ref{expl3}, \ref{expl8}, \ref{expl9}, \ref{expl10} and \ref{expl11}. \Qed
\subsection{The case ${\rm Kod}(S/i)=1$}

Now we give a list of possibilities for the case ${\rm Kod}(S/i)=1.$

\begin{theorem}\label{prop1}
Let $S$ and $P$ be as above.
If ${\rm Kod}(P)=1,$ only the following cases can occur:

\begin{description}
  \item[a)] $\chi(\mathcal{O}_P)=2,$ $q(P)=0$ and
  \begin{description}
    \item[$\cdot$] $\{r_i\}=\emptyset,$ $K_P\overline B=4,$ $\overline B^2=-32,$ $t-8=K_S^2\in\{4,\ldots,8\};$
  \end{description}

  \item[b)] $\chi(\mathcal{O}_P)=1,$ $q(P)=0$ and
  \begin{description}
    \item[$\cdot$] $\{r_i\ne2\}=\emptyset,$ $K_P\overline B=2,$ $\overline B^2=-12,$ $t-2=K_S^2\in\{2,3,4\},$ or
    \item[$\cdot$] $\{r_i\ne2\}=\emptyset,$ $K_P\overline B=4,$ $\overline B^2=-16,$ $t=K_S^2\in\{4,\ldots,8\},$ or
    \item[$\cdot$] $\{r_i\ne2\}=\{4\},$ $K_P\overline B=2,$ $\overline B^2=-4,$ $t=K_S^2\in\{4,\ldots,8\};$
  \end{description}

  \item[c)] $\chi(\mathcal{O}_P)=1,$ $q(P)=1$ and
  \begin{description}
    \item[$\cdot$] $\{r_i\ne2\}=\emptyset,$ $K_P\overline B=2,$ $\overline B^2=-12,$ $t-2=K_S^2\in\{2,\ldots,6\},$ or
    \item[$\cdot$] $\{r_i\}=\emptyset,$ $K_P\overline B=4,$ $\overline B^2=-16,$ $t=K_S^2\in\{4,\ldots,8\}.$
  \end{description}

  \item[d)] $\chi(\mathcal{O}_P)=0,$ $q(P)=1$ and
  \begin{description}
    \item[$\cdot$] $\{r_i\ne2\}=\emptyset,$ $K_P\overline B=2,$ $\overline B^2=4,$ $t=0,$ $K_S^2=6,$ or
    \item[$\cdot$] $\{r_i\ne2\}=\emptyset,$ $K_P\overline B=4,$ $\overline B^2=0,$ $t=0,$ $K_S^2=8,$ or
    \item[$\cdot$] $\{r_i\ne2\}=\{4\},$ $K_P\overline B=2,$ $\overline B^2=12,$ $t=0,$ $K_S^2=8.$
  \end{description}
\end{description}
Furthermore, there exist examples for
\begin{description}
  \item[$\cdot$] {\rm a)} with $K_S^2=8;$
  \item[$\cdot$] {\rm b)} with $K_S^2=4,$ $6$ or $7;$
  \item[$\cdot$] {\rm c)} with $K_S^2=8;$
  \item[$\cdot$] {\rm d)} with $K_S^2=6$ or $8.$
\end{description}
\end{theorem}
{\bf Proof :}
Since $p_g(P),q(P)\leq 1,$ we have the following cases:
\begin{description}
  \item[a)] $\chi(\mathcal O_P)=2,$ $q(P)=0.$\\
  From Proposition \ref{m}, b) it is immediate that
  $h^0(W,\mathcal{O}_W(2K_W+L))=1$ and $K_W^2=0$ (thus $\overline B$ is smooth).
  Proposition \ref{equation} gives $K_P\overline B=4$ and $\overline B^2=-32.$
  If $K_S^2=9,$ then the number of nodal curves of $B$ is $t=K_S^2+8=17,$
  from Proposition \ref{numbernodes}, a).
  This is impossible because Proposition \ref{Miy} implies $t\leq 16.$
  Proposition \ref{numbernodes}, c) gives $K_S^2\geq 4.$
  \item[b)] $\chi(\mathcal O_P)=1,$ $q(P)=0.$\\
  Proposition \ref{m}, a) implies $h^0(W,\mathcal{O}_W(2K_W+L))\leq 3,$
  with equality only if $K_S^2=9$ and $K_W^2=0$
  (hence $\sum(r_i-2)=0$ and $W=P$).
  In this case Proposition \ref{equation}, a) implies
  $K_WB'=6$ and then $B'\ne\emptyset.$
  Now $p_a(B')=1$ (see Claim 1 in the proof of Proposition \ref{m}), thus $B'$ is an union of elliptic
  components. But Proposition \ref{Sa} implies that a minimal surface of general type with $\chi=1$ and
  $K^2=9$ contains no elliptic curves. Therefore $h^0(W,\mathcal{O}_W(2K_W+L))\leq 2.$

  Since ${\rm Kod}(P)=1,$ $K_P\overline B=0$ implies that $\overline B$ is contained
  in the elliptic fibration of $P$ and then $S$ has an elliptic fibration,
  which is impossible because $S$ is of general type.

  So $K_P\overline B\ne 0.$ Now Propositions \ref{equation} and \ref{numbernodes}, a) give
  the following possibilities:
  \begin{description}
   \item[1)] $\sum(r_i-2)=0,$ $K_P\overline B=2,$ $\overline B^2=-12,$ $t=K_S^2+2;$
   \item[2)] $\sum(r_i-2)=0,$ $K_P\overline B=4,$ $\overline B^2=-16,$ $t=K_S^2;$
   \item[3)] $\sum(r_i-2)=2,$ $K_P\overline B=2,$ $\overline B^2=-4,$ $t=K_S^2.$
  \end{description}

  In case 1), $t>6$ implies $\overline{B'}^2=\overline{B}^2+2t>0,$
  a contradiction (see Section \ref{q1}).

  Similarly $t\leq 8,$ in case 2). Proposition \ref{numbernodes}, c) gives $K_S^2\geq 4,$ in this case.

  In case 3), the quadruple or $(3,3)$-point of $\overline B$ gives rise to an elliptic curve
  in $S,$ thus $K_S^2\ne 9,$ from Proposition \ref{Sa}. Again Proposition \ref{numbernodes}, c)
  implies $K_S^2\geq 4.$
  \item[c)] $\chi(\mathcal O_P)=1,$ $q(P)=1.$\\
  This is analogous to the proof of b): just notice that Proposition \ref{m}, b)
  excludes case 3) and implies $K_W^2=0$ in case 2); in case 1) is no longer true that
  $t\leq 6,$ instead use Proposition \ref{Miy} to obtain $t\leq 8$ (thus $K_S^2\leq 6$).
  \item[d)] $\chi(\mathcal O_P)=0,$ $q(P)=1.$\\ As in b), one shows that $h^0(W,\mathcal{O}_W(2K_W+L))\leq 3$
  and $K_P\overline B\ne 0.$ Propositions \ref{equation} and \ref{numbernodes}, a)
  give the following possibilities:
  \begin{description}
   \item[1)] $\sum(r_i-2)=0,$ $K_P\overline B=2,$ $\overline B^2=4,$ $t=K_S^2-6;$
   \item[2)] $\sum(r_i-2)=0,$ $K_P\overline B=4,$ $\overline B^2=0,$ $t=K_S^2-8;$
   \item[3)] $\sum(r_i-2)=2,$ $K_P\overline B=2,$ $\overline B^2=12,$ $t=K_S^2-8.$
  \end{description}

  As in the proof of b), the existence of a quadruple or $(3,3)$-point on $\overline B$ implies
  $K_S^2\ne 9,$ in case 3).

  Consider now cases 1) and 2). From Proposition \ref{Miy}, $P$ has no smooth rational curves.
  Any singular rational curve $D$ of $\overline B$ satisfies $D^2\leq 0,$ because, since $\overline B$
  has only negligible singularities, $D$ is contained in fibres of a fibration $\overline{f_A}$ of $P$
  (see Section \ref{q1}).
  Therefore $B$ has no $(-2)$-curves, {\it i.e.} $t=0.$
\end{description}
The examples are given in Sections \ref{expl4}, \ref{expl5}, \ref{expl6}, \ref{expl9},
\ref{expl10} and \ref{expl7}. \Qed

\subsection{The case ${\rm Kod}(S/i)=2$}

Finally we give a list of possibilities for the case ${\rm Kod}(S/i)=2.$

\begin{theorem}\label{prp2}
Let $S$ and $P$ be as above.
If ${\rm Kod}(P)=2,$ then $\overline B$ has at most negligible singularities and only the following cases can occur:

\begin{description}
  \item[a)] $\chi(\mathcal{O}_P)=2,$ $q(P)=0$ and
  \begin{description}
    \item[$\cdot$] $K_P\overline B=0,$ $\overline B^2=-24,$
    $t=12,$ $K_S^2=2K_P^2,$ $K_P^2=2,3,4,$ or
    \item[$\cdot$] $K_P\overline B=2,$ $\overline B^2=-28,$
    $t-10+2K_P^2=K_S^2\in\{2K_P^2+2,\ldots,2K_P^2+4\},$ $K_P^2=1,2;$
  \end{description}

  \item[b)] $\chi(\mathcal{O}_P)=1,$ $q(P)=1$ and
  \begin{description}
    \item[$\cdot$] $K_P\overline B=0,$ $\overline B^2=-8,$
    $t=4,$ $K_S^2=2K_P^2,$ $K_P^2=2,3,4,$ or
    \item[$\cdot$] $K_P\overline B=2,$ $\overline B^2=-12,$
    $K_P^2=2,$ $t+2=K_S^2\in\{6,7,8\};$
  \end{description}

  \item[c)] $\chi(\mathcal{O}_P)=1,$ $q(P)=0$ and
  \begin{description}
    \item[$\cdot$] $K_P\overline B=0,$ $\overline B^2=-8,$
    $t=4,$ $K_S^2=2K_P^2,$ $K_P^2=1,\ldots,4,$ or
    \item[$\cdot$] $K_P\overline B=2,$ $\overline B^2=-12,$
    $t+2K_P^2-2=K_S^2\in\{2K_P^2+2,\ldots,2K_P^2+4\},$ $K_P^2=1,2,$ or
    \item[$\cdot$] $K_P\overline B=4,$ $\overline B^2=-16,$
    $K_P^2=1,$ $t+2=K_S^2\in\{6,7,8\}.$
  \end{description}

\end{description}
Moreover, there exist examples for
\begin{description}
  \item[$\cdot$] {\rm a)} with $K_S^2=4,$ $6,$ $7$ or $8;$
  \item[$\cdot$] {\rm b)} with $K_S^2=4.$
\end{description}
\end{theorem}
{\bf Proof :} \\\\
{\bf Claim :} If $K_P\overline B=0,$ then $\overline B$ is a disjoint union
of nodal curves.\\\\
{\em Proof}\ : As $P$ is of general type, $\overline B$ is an union of
nodal curves. Suppose that $\overline B$ has a singularity. Then it contains two
nodal curves $D_1,$ $D_2$ such that $D_1D_2\geq 2,$ otherwise $B$ has a $(-3)$-curve,
contradicting $B\equiv 0\ ({\rm mod\ }2).$ Since $K_P^2>0,$ $K_P(D_1+D_2)=0$ and
$(D_1+D_2)^2\not <0,$ the index theorem implies that $D_1+D_2$ is homologous to zero,
a contradiction.\ $\diamondsuit$
\\

Since $p_g(P),q(P)\leq 1$ and $p_g(P)\geq q(P),$ we have only the following
three cases:
\begin{description}
  \item[a)] $\chi(\mathcal O_P)=2,$ $q(P)=0.$\\
  Propositions \ref{equation}, a) and \ref{m}, b) give:
  $$h^0(W,\mathcal O_W(2K_W+L))=K_P^2+K_P\delta+\frac{1}{2}\sum(r_i-2)-1,$$
  $$h^0(W,\mathcal O_W(2K_W+L))\leq \frac{1}{2}K_W^2+1\leq \frac{1}{2}K_P^2+1.$$
  From this we get
  \begin{equation}\label{idem}
  \frac{1}{2}K_P^2+K_P\delta+\frac{1}{2}\sum(r_i-2)\leq 2,
  \end{equation}
  with equality only if $K_W^2=K_P^2.$
  Since $K_P^2>0,$ $K_P\delta=0$ or $1.$

  If $K_P\delta=1,$ then $\sum(r_i-2)=0$ and $K_P^2=1$ or $2.$

  If $K_P\delta=0,$ then $\sum(r_i-2)=0,$ from the Claim above.
  As $K_S^2\leq 9,$ Proposition \ref{numbernodes} implies $h^0(W,\mathcal O_W(2K_W+L))\leq 3.$

  Now the result follows from Propositions \ref{equation} and \ref{numbernodes}, a).
  Notice that Proposition \ref{equation} gives \ $\overline B^2\geq -2(12+2K_P\delta).$
  This implies {$t\leq 12+2K_P\delta,$} because, since $q(P)=0$
  and $\overline B$ has only negligible singularities, every component of $\overline B$
  has non-positive self-intersection.

  \item[b)] $\chi(\mathcal O_P)=1,$ $q(P)=1.$\\
  Equation (\ref{idem}) is still valid here.
  As $p_g(P)=q(P)=1,$ then $K_P^2\geq 2,$ hence $$K_P\delta+\frac{1}{2}\sum(r_i-2)\leq 1.$$
  Using the Claim above and Proposition \ref{equation}, we have\\
  $K_P\delta=0,\ \sum(r_i-2)=0,\ K_P^2=2,3$ or $4,$ $\overline B^2=-8,\ t=4,$ or\\
  $K_P\delta=1,\ \sum(r_i-2)=0,\ K_P^2=2,$ $\overline B^2=-12.$\\
  Now the result follows from Proposition \ref{numbernodes}, a)
  (notice that $K_P^2=2$ implies $t\ne 7,$ by Theorem \ref{Miy}).

  \item[c)] $\chi(\mathcal O_P)=1,$ $q(P)=0.$\\
  Propositions \ref{equation}, a), \ref{numbernodes}, c) and \ref{m}, a) imply
  \begin{equation}\label{a1}
  h^0(W,\mathcal O_W(2K_W+L))=K_P(K_P+\delta)+\frac{1}{2}\sum(r_i-2)\leq 4
  \end{equation}
  and
  \begin{equation}\label{a2}
  h^0(W,\mathcal O_W(2K_W+L))\leq 3+\frac{1}{3}K_W^2,
  \end{equation}
  with equality only if $K_S^2=9.$ As $K_P^2\geq 1$ and $K_W^2\leq K_P^2,$
  this implies $K_P\delta\leq 2.$

  $\bullet$ Suppose $K_P\delta=0.$

  We have $\sum(r_i-2)=0,$ by the Claim above. Hence
  $h^0(W,\mathcal O_W(2K_W+L))=K_P^2\leq 4,$ by (\ref{a1}). Now from Proposition
  \ref{equation}, b) and Proposition \ref{numbernodes}, b), we have
  $\overline B^2=(2\delta)^2=-8$ and $t\geq 4.$ Thus $t=4$ and, using Proposition
  \ref{numbernodes}, a), we conclude that
  $$K_S^2=2K_P^2,\ 1\leq K_P^2\leq 4.$$

  $\bullet$ Suppose $K_P\delta=1.$

  Then $h^0(W,\mathcal O_W(2K_W+L))=4$ only if $K_S^2=9,$
  from (\ref{a1}) and (\ref{a2}). 
  In this case Propositions \ref{numbernodes}, b) and \ref{equation}, a) imply
  $K_P^2=3$ and $B$ contains an elliptic curve, which contradicts Proposition \ref{Sa}.

  If  $K_P^2=1/2\sum(r_i-2)=1,$ then $K_W^2=0$ and $K_S^2=9,$
  by (\ref{a2}). The quadruple or $(3,3)$-point of $\overline B$ gives rise to an elliptic curve in $S,$
  which is impossible from Proposition \ref{Sa}.

  Now using (\ref{a1}), Proposition
  \ref{equation}, b) and Proposition \ref{numbernodes}, we obtain
  $$K_P^2=2,\ \sum(r_i-2)=0,\ \delta^2=-3,\ t=K_S^2-2\geq 4$$
  or
  $$K_P^2=1,\ \sum(r_i-2)=0,\ \delta^2=-3,\ t=K_S^2\geq 4.$$
  Hence $\overline B^2=-12$ and then $t\leq 6,$ because, since $q(P)=0$
  and $\overline B$ has only negligible singularities, every component of $\overline B$
  has non-positive self-intersection.

  $\bullet$ Suppose $K_P\delta=2.$

  Then $K_P^2\leq 2,$ from (\ref{a1}), and
  $h^0(W,\mathcal O_W(2K_W+L))\leq 3,$ from (\ref{a2}).
  The only possibility allowed by (\ref{a1}), Proposition \ref{equation}, b)
  and Proposition \ref{numbernodes} is:
  $$K_P^2=1,\ \delta^2=-4,\ t=K_S^2-2\geq 4.$$

  It remains to be shown that $K_S^2\ne 9.$ In this case, the curve $\overline B$
  has at least 8 disjoint components contained in a fibration $\overline{f_A}$ of $P$
  (see Section \ref{q1}). These are independent in ${\rm Pic}(P)$ from a general
  fibre of $\overline{f_A}$ and from $K_P,$ so ${\rm Pic}(P)$ has 10
  independent classes. This is a contradiction because the second Betti number of $P$ is
  $$b_2(P)=12\chi(\mathcal O_P)-K_P^2+4q(P)-2=9.$$
\end{description}
The examples can be found in Sections \ref{expl5}, \ref{expl3}, \ref{expl6},
\ref{expl8}, \ref{expl9}, \ref{expl10} and \ref{expl11}. \Qed

\section{(Bi)double cover examples}\label{bdce}

The next sections contain constructions of minimal smooth surfaces of general type $S$ with
$p_g=q=1$ which present examples for cases a), b), ...
of Theorems \ref{prop0}, \ref{prop1} and \ref{prp2}.
Only the existence of case c) of Theorem \ref{prp2} remains an open problem.

Each example is obtained as the smoth minimal model of a bidouble cover of a ruled
surface (irregular only in Sections \ref{expl7} and \ref{expl11}).
\\

A bidouble cover is a finite flat Galois morphism with Galois group $\mathbb Z_2^2.$
Following \cite{Ca2} or \cite{Pa}, to define a bidouble cover cover $\psi:V\rightarrow X,$
with $V,$ $X$ smooth surfaces, it suffices to present:
\begin{description}
\item[$\cdot$] smooth divisors $D_1, D_2, D_3\subset X$ with pairwise
transverse intersections and no common intersection;
\item[$\cdot$] line bundles $L_1, L_2, L_3$ such that $2L_g\equiv D_j+D_k$ for
each permutation $(g,j,k)$ of $(1,2,3).$
\end{description}
If ${\rm Pic}(X)$ has no 2-torsion, the $L_i$'s are uniquely determined
by the $D_i$'s.

Let $N:=2K_X+\sum_1^3L_i.$ One has:
$$p_g(V)=p_g(X)+\sum_1^3h^0(X,\mathcal O_X(K_X+L_i)),$$
$$\chi(\mathcal O_V)=4\chi(\mathcal O_X)+\frac{1}{2}\sum_1^3 L_i(K_X+L_i),$$
$$2K_V\equiv\psi^*\left(N\right)$$
and
$$H^0(V,\mathcal O_V(2K_V))\simeq H^0(X,\mathcal O_X(N))\oplus
\bigoplus_{i=1}^3 H^0(X,\mathcal O_X(N-L_i)).$$

The bicanonical map of $V$ is composed with the involution $i_g,$ associated to $L_g,$
if and only if $$h^0(X,\mathcal O_X(2K_X+L_g+L_j))=h^0(X,\mathcal O_X(2K_X+L_g+L_k))=0.$$

For more information on bidouble covers see \cite{Ca2} or \cite{Pa}.\\

Denote by $i_1,i_2,i_3$ the involutions of $V$ corresponding to $L_1,L_2,L_3,$ respectively.
In each example the invariants of $W_{j}:=V/{i_j},$ $j=1,2,3,$ are given,
but the details are not written here. These can be found in the author's Ph.D. thesis,
which can be sent by e-mail under request.\\

We use the following:

\begin{notation}\label{notation}
Let $p_0,\ldots,p_j,\ldots,p_{j+s}\in\mathbb P^2$ be distinct points and define
$T_i$ as the line through $p_0$ and $p_i,$ $i=1,\ldots,j.$
We say that a plane curve is of type $$d(m,(n,n)_T^j,r^s)$$
if it is of degree $d$ and if it has: an $m$-uple point at $p_0,$
an $(n,n)$-point at $p_1,\ldots,p_j,$ an $r$-uple point at $p_{j+1},\ldots,p_{j+s}$
and no other non-negligible singularities.
The index $_T$ is used if $T_i$ is tangent to the $(n,n)$-point at $p_i.$

An obvious generalization is used if there are other singularities.\\

Let $p_1',\ldots,p_j'$ be the infinitely near points to $p_1,\ldots,p_j,$ respectively.
We denote by $$\mu:X\rightarrow\mathbb P^2$$ the blow-up with centers
$$p_0,p_1,p_1',\ldots,p_j,p_j',p_{j+1},\ldots,p_{j+s}$$ and by
$$E_0,E_1,E_1',\ldots,E_j,E_j',E_{j+1},\ldots,E_{j+s}$$
the corresponding exceptional divisors (with self-intersection $-1$).

The notation $\widetilde{\cdot}$ stands for the total transform $\mu^*(\cdot)$ of a curve.

The letter $T$ is reserved for a general line of $\mathbb P^2.$
\end{notation}
The genus of a general Albanese fibre of $S$ is denoted by $g.$
\subsection[$K^2=8,$ $g=3,$ $S/{i_1}$ ruled, $S/{i_2}$ rational, ${\rm Kod}(S/{i_3})=1$]
{$K^2=8,$ $g=3,$\\
$S/{i_1}$ ruled, $S/{i_2}$ rational, ${\rm Kod}(S/{i_1})=1$}\label{expl4}

Here we construct a surface of general type $V$ with $p_g=q=1$ and $g=3$ such that
$K_S^2=8,$ where $S$ is the minimal model of $V.$
The quotients $W_j:=V/{i_j},$ $j=1,2,3,$ satisfy:
\begin{description}
\item[$\cdot$] $W_1$ is ruled, $q(W_1)=1;$
\item[$\cdot$] $W_2$ is rational;
\item[$\cdot$] ${\rm Kod}(W_3)=p_g(W_3)=1,$ $q(W_3)=0$
\end{description}

This gives an example for case a) of Theorem \ref{prop1}.
The surface $S$ is a Du Val double plane of type III described in \cite{Po3}.\\\\
{\bf$\cdot$ Construction of $S$}
\\
Let $Q$ be a reduced curve of type $4(0,(2,2)_T^2),$ {\it i.e.} $Q$ is the union of
two conics tangent to the lines $T_1$ and $T_2$ at $p_1,$ $p_2.$
Let $C$ be another non-degenerate conic tangent to $T_1,$ $T_2$ at
$p_1,$ $p_2$ and let $T_3,\ldots,T_6\ne T_1, T_2$ be distinct lines
through $p_0\in T_1\bigcap T_2.$

Set:

$$
\begin{array}{l}
D_1:=\widetilde{T_1}+\cdots+\widetilde{T_6}-\sum_1^2(E_i+E_i')-6E_0,\\
D_2:=\widetilde{Q}-\sum_1^2(E_i+3E_i'),\\
D_3:=\widetilde{C}-\sum_1^2(E_i+E_i').
\end{array}
$$

We have
$$
\begin{array}{l}
L_1\equiv 3\widetilde T-\sum_1^2(E_i+2E_i'),\\
L_2\equiv 4\widetilde T-\sum_1^2(E_i+E_i')-3E_0,\\
L_3\equiv 5\widetilde T-\sum_1^2(E_i+2E_i')-3E_0.
\end{array}
$$
and
$$
\begin{array}{l}
K_X+L_1\equiv -E_1'-E_2'+E_0,\\
K_X+L_2\equiv \widetilde T-2E_0,\\
K_X+L_3\equiv 2\widetilde T-E_1'-E_2'-2E_0.
\end{array}
$$
Let $V\rightarrow X$ be the bidouble cover determined by $D_1,$ $D_2,$ $D_3$
and $S$ be the minimal model of $V.$ We have
$$p_g(S)=\sum_1^3h^0(X,\mathcal O_X(K_X+L))=0+0+1=1,$$
$$\chi(\mathcal O_S)=4+\frac{1}{2}\sum_1^3 L_i(K_X+L_i)=4-2-1+0=1.$$

One has
$$\begin{array}{l}
h^0(X,\mathcal O_X(2K_X+L_1+L_2))=0,\\
h^0(X,\mathcal O_X(2K_X+L_1+L_3))=1,\\
h^0(X,\mathcal O_X(2K_X+L_2+L_3))=0,
\end{array}
$$
thus the bicanonical map of $V$ is composed with the involution $i_2$ and is not
composed with the involutions $i_1$ and $i_3.$

The {\bf Albanese fibration} of $S$
is given by the pullback of the pencil of conics tangent to $T_1,T_2$ at $p_1,p_2.$
\\\\
{\bf$\cdot$ Calculation of $K_S^2$}
\\
We have
$$N:=2K_X+\sum_1^3L_i\equiv N_1+N_2,$$
where $$N_1:=\widetilde{T_1}+\widetilde{T_2}-2E_0-\sum_1^22E_i',\hspace{2em} N_2:=4\widetilde T-2E_0-\sum_1^2(E_i+E_i').$$
Since the support of the pullback of $N_1$ is a disjoint union of $8=-N_1^2$ $(-1)$-curves,
$|N_2|$ has no fixed component and $N_1N_2=0,$
$$K_S^2=K_V^2+(-N_1^2)=N^2-N_1^2=N_2^2=8.$$

\subsection[$K^2=6,$ $g=4,$ ${\rm Kod}(S/{i_1})=2,$ $S/{i_2}$ rational, ${\rm Kod}(S/{i_3})=1$]
{$K^2=6,$ $g=4,$\\
${\rm Kod}(S/{i_1})=2,$ $S/{i_2}$ rational, ${\rm Kod}(S/{i_3})=1$}\label{expl5}

This section contains the construction of a bidouble cover
$V\rightarrow X$ such that the minimal model $S$ of $V$ is a
surface of general type with $K^2=6,$ $p_g=q=1,$
$g=4$ and that the quotients $W_j:=V/{i_j}$ satisfy:
\begin{description}
\item[$\cdot$] ${\rm Kod}(W_1)=2,$ $p_g(W_1)=1,$ $q(W_1)=0;$
\item[$\cdot$] $W_2$ is rational;
\item[$\cdot$] ${\rm Kod}(W_3)=1,$ $p_g(W_3)=0,$ $q(W_3)=1.$
\end{description}

This is an example for Theorems \ref{prop1}, d) and \ref{prp2}, a).

One can verify that $S$ is a Du Val double plane obtained imposing a $4$-uple
point to the branch locus of a Du Val's ancestor of type $\mathcal D_5$
({\it cf.} \cite{Ri2}).

Recall Notation \ref{notation}.
\\\\
{\bf$\cdot$ Construction of $S$}
\\
From Proposition \ref{PropUsefulPencils} in Appendix \ref{UsefulPencils},
there is a pencil $l,$ with no base component, of curves of type
$7(3,(2,2)_T^5).$ Let $Q$ be a general element of this pencil and $C$
be a reduced curve of type $4(2,(1,1)_T^5).$

Set:

$$
\begin{array}{l}
D_1:=\widetilde{T_1}+\cdots+\widetilde{T_4}-\sum_1^4(E_i+E_i')+(E_5-E_5')-4E_0,\\
D_2:=\widetilde{T_5}+\widetilde{Q}-\sum_1^4(E_i+3E_i')-3E_5-3E_5'-4E_0,\\
D_3:=\widetilde{C}-\sum_1^5(E_i+E_i')-2E_0.
\end{array}
$$

We have
$$
\begin{array}{l}
L_1\equiv 6\widetilde T-\sum_1^4(E_i+2E_i')-2E_5-2E_5'-3E_0,\\
L_2\equiv 4\widetilde T-\sum_1^4(E_i+E_i')-E_5'-3E_0,\\
L_3\equiv 6\widetilde T-\sum_1^5(E_i+2E_i')-4E_0
\end{array}
$$
and
$$
\begin{array}{l}
K_X+L_1\equiv \left(2\widetilde T-\sum_1^4E_i'-E_0\right)+ \left(\widetilde{T_5}-E_5-E_5'-E_0\right),\\
K_X+L_2\equiv \widetilde T+E_5-2E_0,\\
K_X+L_3\equiv 3\widetilde T-\sum_1^5E_i'-3E_0,
\end{array}
$$
hence
$$p_g(S)=\sum_1^3h^0(X,\mathcal O_X(K_X+L_i))=1+0+0=1$$ and
$$\chi(\mathcal O_S)=4+\frac{1}{2}\sum_1^3 L_i(K_X+L_i)=4+0-1-2=1.$$

One has
$$\begin{array}{l}
h^0(X,\mathcal O_X(2K_X+L_1+L_2))=0,\\
h^0(X,\mathcal O_X(2K_X+L_1+L_3))=2,\\
h^0(X,\mathcal O_X(2K_X+L_2+L_3))=0,
\end{array}
$$
thus the bicanonical map of $V$ is composed with the involution $i_2$ and is not
composed with the involutions $i_1$ and $i_3.$

The {\bf Albanese fibration} of $S$
is given by the pullback of the pencil $l.$
\\\\
{\bf$\cdot$ Calculation of $K_S^2$}
\\
We have
$$N:=2K_X+\sum_1^3L_i\equiv N_1+N_2,$$
where $$N_1:=\widetilde{T_1}+\cdots+\widetilde{T_5}-5E_0-\sum_1^52E_i',\ \ \ \ N_2:=5\widetilde T-3E_0-\sum_1^5(E_i+E_i').$$
As the support of the pullback of $N_1$ is a disjoint union of $20=-N_1^2$ $(-1)$-curves,
$|N_2|$ has no fixed component and $N_1N_2=0,$ then
$$K_S^2=K_V^2+(-N_1^2)=N^2-N_1^2=N_2^2=6.$$

\subsection[$K^2=4,$ $g=3,$ ${\rm Kod}(S/{i_1})=2,$ $S/{i_2}$ rational, ${\rm Kod}(S/{i_3})=0$]
{$K^2=4,$ $g=3,$\\
 ${\rm Kod}(S/{i_1})=2,$ $S/{i_2}$ rational, ${\rm Kod}(S/{i_3})=0$}\label{expl3}

In this section a surface of general type $S$ with $p_g=q=1,$ $K^2=4$ and
$g=3$ is constructed. It is the minimal model of a double cover of surfaces
$W_1,$ $W_2,$ $W_3$ such that:
\begin{description}
\item[$\cdot$] ${\rm Kod}(W_1)=2,$ $p_g(W_1)=1,$ $q(W_1)=0;$
\item[$\cdot$] $W_2$ is rational;
\item[$\cdot$] ${\rm Kod}(W_3)=0,$ $p_g(W_3)=0,$ $q(W_3)=1.$
\end{description}

This gives an example for Theorems \ref{prop0}, b) and \ref{prp2}, a).

One can verify that $S$ is a Du Val double plane obtained imposing two $4$-uple
points to the branch locus of a Du Val's ancestor of type $\mathcal D_4$
({\it cf.} \cite{Ri2}).

We keep Notation \ref{notation}.
\\\\
{\bf$\cdot$ Construction of $S$}
\\
From Proposition \ref{PropUsefulPencils},
there is a pencil $l,$ with no base component, of curves of type
$6(2,(2,2)_T^4),$ through points $p_0,\ldots,p_4$ 
({\it i.e.} of plane curves of degree $6$ with a double point at $p_0$ and
a tacnode at $p_i$ with tangent line through $p_0,p_i,$ $i=1,\ldots,4$). 
Let $Q$ be a general element of this pencil, $C$ be a
reduced curve of type $4(2,(1,1)_T^4)$ and set:
$$
\begin{array}{l}
D_1:=\widetilde{T_1}+\cdots+\widetilde{T_4}-\sum_1^4(E_i+E_i')-4E_0,\\
D_2:=\widetilde Q-\sum_1^4(E_i+3E_i')-2E_0,\\
D_3:=\widetilde C-\sum_1^4(E_i+E_i')-2E_0.
\end{array}
$$
Let $V\rightarrow X$ be the bidouble cover determined by $D_1,$ $D_2,$ $D_3$
and $S$ be the minimal model of $V.$

We have
$$
\begin{array}{l}
L_1\equiv 5\widetilde T-\sum_1^4(E_i+2E_i')-2E_0,\\
L_2\equiv 4\widetilde T-\sum_1^4(E_i+E_i')-3E_0,\\
L_3\equiv 5\widetilde T-\sum_1^4(E_i+2E_i')-3E_0
\end{array}
$$
and
$$
\begin{array}{l}
K_X+L_1\equiv 2\widetilde T-\sum_1^4E_i'-E_0,\\
K_X+L_2\equiv \widetilde T-2E_0,\\
K_X+L_3\equiv 2\widetilde T-\sum_1^4E_i'-2E_0.
\end{array}
$$
Then
$$p_g(S)=\sum_1^3h^0(X,\mathcal O_X(K_X+L_i))=1+0+0=1$$ and
$$\chi(\mathcal O_S)=4+\frac{1}{2}\sum_1^3 L_i(K_X+L_i)=4+0-1-2=1.$$

As
$$\begin{array}{l}
h^0(X,\mathcal O_X(2K_X+L_1+L_2))=0,\\
h^0(X,\mathcal O_X(2K_X+L_1+L_3))=1,\\
h^0(X,\mathcal O_X(2K_X+L_2+L_3))=0,
\end{array}
$$
the bicanonical map of $V$ is composed with the involution $i_2,$ corresponding
to $L_2,$ and is not composed with the involutions $i_1,i_3,$ corresponding
to $L_1,L_3.$

The {\bf Albanese fibration} of $S$
is given by the pullback of the pencil $l.$
\\\\
{\bf$\cdot$ Calculation of $K_S^2$}
\\
We have
$$N:=2K_X+\sum_1^3L_i\equiv N_1+N_2,$$
where $$N_1:=\widetilde{T_1}+\cdots+\widetilde{T_4}-4E_0-\sum_1^42E_i',\ \ \ \ N_2:=4\widetilde T-2E_0-\sum_1^4(E_i+E_i').$$
As the support of the pullback of $N_1$ is a disjoint union of $16=-N_1^2$ $(-1)$-curves,
$|N_2|$ has no fixed component and $N_1N_2=0,$
$$K_S^2=K_V^2+(-N_1^2)=N^2-N_1^2=N_2^2=4.$$

\subsection[$K^2=4,$ $g=2,$ $S/{i_1}$ ruled, ${\rm Kod}(S/{i_2})=1,$ ${\rm Kod}(S/{i_3})=2$]
{$K^2=4,$ $g=2,$\\
$S/{i_1}$ ruled, ${\rm Kod}(S/{i_2})=1,$ ${\rm Kod}(S/{i_3})=2$}\label{expl6}

This section contains the construction of a surface of
general type $S$ with $p_g=q=1,$ $K^2=4,$ $g=2$ and ${\rm deg}(\phi_2)=2,$
where $\phi_2$ is the bicanonical map of $S.$
The surface $S$ is the minimal model of a double cover
of surfaces $W_1,$ $W_2,$ $W_3$ such that:
\begin{description}
\item[$\cdot$] $W_1$ is ruled, $q(W_1)=1;$
\item[$\cdot$] ${\rm Kod}(W_2)=1,$ $p_g(W_2)=q(W_2)=0;$
\item[$\cdot$] ${\rm Kod}(W_3)=2,$ $p_g(W_3)=1,$ $q(W_3)=0.$
\end{description}

This gives an example for Theorems \ref{prop1}, b) and \ref{prp2}, a).

Notation \ref{notation} is used in this section.
\\\\
{\bf$\cdot$ Construction of $S$}
\\
By Proposition \ref{PropUsefulPencils} in Appendix \ref{UsefulPencils},
there is a pencil $l,$ with no base component, of curves of type
$6(2,(2,2)_T^4).$ Let $Q_1$ be a general element of this pencil, $Q_2$ be a smooth curve
of type $3(1,(1,1)_T^4)$ and $Q:=Q_1+Q_2.$ Let $T_5$ be a line through $p_0$
transverse to $Q.$

Set:
$$
\begin{array}{l}
D_1:=\widetilde{T_1}+\widetilde Q-4E_1-4E_1'-\sum_2^4(3E_i+3E_i')-4E_0,\\
D_2:=\widetilde{T_2}+\cdots+\widetilde{T_5}-\sum_2^4(E_i+E_i')-4E_0,\\
D_3:=\sum_2^4(E_i-E_i').
\end{array}
$$
We have
$$
\begin{array}{l}
L_1\equiv 2\widetilde T-\sum_2^4 E_i'-2E_0,\\
L_2\equiv 5\widetilde T-2E_1-2E_1'-\sum_2^4(E_i+2E_i')-2E_0,\\
L_3\equiv 7\widetilde T-\sum_1^4(2E_i+2E_i')-4E_0
\end{array}
$$
and
$$
\begin{array}{l}
K_X+L_1\equiv -\widetilde T+\sum_1^4E_i+E_1'-E_0,\\
K_X+L_2\equiv 2\widetilde T-\sum_1^4E_i'-E_1-E_0,\\
K_X+L_3\equiv \widetilde{T_1}+\cdots+\widetilde{T_4}-\sum_1^4(E_i+E_i')-3E_0.
\end{array}
$$
Let $\psi:V\rightarrow X$ be the bidouble cover determined by $D_1,$ $D_2,$ $D_3$
and $S$ be the minimal model of $V.$ Then
$$p_g(S)=\sum_1^3h^0(X,\mathcal O_X(K_X+L_i))=0+0+1=1$$ and
$$\chi(\mathcal O_S)=4+\frac{1}{2}\sum_1^3 L_i(K_X+L_i)=4-2-1+0=1.$$
One can verify that
$$\begin{array}{l}
h^0(X,\mathcal O_X(2K_X+L_1+L_2))=0,\\
h^0(X,\mathcal O_X(2K_X+L_1+L_3))=0,\\
h^0(X,\mathcal O_X(2K_X+L_2+L_3))=1,
\end{array}
$$
thus the bicanonical map of $V$ is composed with the involution $i_1$ and is not
composed with the involutions $i_2$ and $i_3.$

The {\bf Albanese fibration} of $S$
is given by the pullback of the pencil of lines through $p_0.$
\\\\
{\bf$\cdot$ Calculation of $K_S^2$}
\\
We have
$$N:=2K_X+\sum_1^3L_i\equiv N_1+N_2,$$
where $$N_1:=\widetilde{T_1}+\cdots+\widetilde{T_4}-4E_0-(E_1+E_1')-\sum_2^4 2E_i'$$
and $$N_2:=4\widetilde T-2E_0-\sum_1^4(E_i+E_i').$$
As the support of the pullback of $N_1$ is a disjoint union of $14=-N_1^2$ $(-1)$-curves,
$|N_2|$ has no fixed component and $N_1N_2=0,$ then
$$K_S^2=N^2+(-N_1^2)=N_2^2=4.$$
\\\\
{\bf$\cdot$ Degree of $\phi_2$}
\\
The system $|\psi^*(N)|$ is strictly contained in the bicanonical system of $V.$
Since $\phi_2$ is composed with $i_1$ and the map $\tau:X\dashrightarrow\mathbb P^2$ induced by
$|N|$ is birational (this can be verified using the Magma function {\em IsInvertible}),
one has ${\rm deg}(\phi_2)=2.$
\subsection[$K^2=8,$ $g=3,$ ${\rm Kod}(S/{i_1})=2,$ ${\rm Kod}(S/{i_2})=0,$ ${\rm Kod}(S/{i_3})=0$]
{$K^2=8,$ $g=3,$\\
${\rm Kod}(S/{i_1})=2,$ ${\rm Kod}(S/{i_2})=0,$ ${\rm Kod}(S/{i_3})=0$}\label{expl8}

A smooth projective surface S of general type is said to be a {\em standard isotrivial
fibration} if there exists a finite group $G$ which acts faithfully on two smooth projective curves
$C$ and $F$ so that $S$ is isomorphic to the minimal desingularization of $T := (C\times F )/G.$
The paper \cite{Po1} contains examples of such surfaces with $K^2=8.$

This section contains the construction of the first surface of general type $S$
with $p_g=q=1,$ $K^2=8$ and $g=3$ which is not a standard isotrivial fibration.

The surface $S$ is the minimal model of a double cover
of surfaces $W_1,$ $W_2,$ $W_3$ such that:
\begin{description}
\item[$\cdot$] ${\rm Kod}(W_1)=2,$ $p_g(W_1)=1,$ $q(W_1)=0;$
\item[$\cdot$] ${\rm Kod}(W_2)=0,$ $p_g(W_2)=0,$ $q(W_2)=1;$
\item[$\cdot$] ${\rm Kod}(W_3)=0,$ $p_g(W_3)=0,$ $q(W_3)=0.$
\end{description}

This is an example for Theorems \ref{prop0} a), b) and \ref{prp2} a).

We use Notation \ref{notation}.
\\\\
{\bf$\cdot$ Construction of $S$}
\\
Let $G$ be a curve of type $6(2,(2,2)_T^4)$ and $C$ be a curve
of type $8(4,(2,2)_T^4,(3,3))$ such that $G+C$ is reduced and
the $(3,3)$-point of $C$ is tangent to $G.$
The existence of these curves is shown in Appendix \ref{CuM}.

Set:
$$
\begin{array}{l}
D_1:=\widetilde{T_1}+\widetilde{T_2}-\sum_1^2 2E_i'+(E_5-E_5')-2E_0,\\
D_2:=\widetilde{G}-\sum_1^4(2E_i+2E_i')-(E_5+E_5')-2E_0,\\
D_3:=\widetilde{T_3}+\widetilde{T_4}+\widetilde{C}-\sum_1^2(2E_i+2E_i')-\sum_3^4(2E_i+4E_i')-(3E_5+3E_5')-6E_0.
\end{array}
$$
We have
$$
\begin{array}{l}
L_1\equiv 8\widetilde T-\sum_1^2(2E_i+2E_i')-\sum_3^4(2E_4+3E_4')-(2E_5+2E_5')-4E_0,\\
L_2\equiv 6\widetilde T-\sum_1^5(E_i+2E_i')-4E_0,\\
L_3\equiv 4\widetilde T-\sum_1^2(E_i+2E_i')-\sum_3^4(E_i+E_i')-E_5'-2E_0
\end{array}
$$
and
$$
\begin{array}{l}
K_X+L_1\equiv 5\widetilde T-\sum_1^2(E_i+E_i')-\sum_3^4(E_4+2E_4')-(E_5+E_5')-3E_0,\\
K_X+L_2\equiv 3\widetilde T-\sum_1^5E_i'-3E_0,\\
K_X+L_3\equiv \widetilde T-\sum_1^2E_i'+E_5-E_0.
\end{array}
$$
Let $V\rightarrow X$ be the bidouble cover determined by $D_1,$ $D_2,$ $D_3$
and $S$ be the minimal model of $V.$ Then
$$\chi(\mathcal O_S)=4+\frac{1}{2}\sum_1^3 L_i(K_X+L_i)=4+0-2-1=1.$$
The Magma procedure $LinSys$ defined in \cite{Ri2} can be used to confirm that
$$h^0(X,\mathcal O_X(K_X+L_1))=1$$ and
$$\begin{array}{l}
h^0(X,\mathcal O_X(2K_X+L_1+L_2))=0,\\
h^0(X,\mathcal O_X(2K_X+L_1+L_3))=1,\\
h^0(X,\mathcal O_X(2K_X+L_2+L_3))=2,
\end{array}
$$
therefore
$$p_g(S)=\sum_1^3h^0(X,\mathcal O_X(K_X+L_i))=1+0+0=1$$ and
the bicanonical map of $V$ is not composed with any of the involutions $i_1,$ $i_2,$ $i_3.$

The {\bf Albanese fibration} of $S$ is induced by a pencil of curves of type\\
$14(6,(4,4)^4_T,(4,4)),$ which contains a fibre equal to $G+C$
(see Appendix \ref{CuM}). From \cite[Theorem 3.2]{Po2}, the existence
of such reducible fibre implies that $S$ is not a standard isotrivial fibration,
so this is not one of Polizzi's examples.
\\\\
{\bf$\cdot$ Calculation of $K_S^2$}
\\
We have
$$N:=2K_X+\sum_1^3L_i\equiv N_1+N_2,$$
where $$N_1:=\widetilde{T_1}+\cdots+\widetilde{T_4}-4E_0-\sum_1^4 2E_i'+(E_5-E_5')$$
and $$N_2:=8\widetilde T-4E_0-\sum_1^5(2E_i+2E_i').$$
As the support of the pullback of $N_1$ is a disjoint union of $18=-N_1^2$ $(-1)$-curves,
$|N_2|$ has no fixed component (again this can be confirmed with the Magma procedure $LinSys$)
and $N_1N_2=0,$ $$K_S^2=N^2+(-N_1^2)=N_2^2=8.$$

\subsection[$K^2=7,$ $g=3,$ ${\rm Kod}(S/{i_1})=2,$ ${\rm Kod}(S/{i_2})=1,$ ${\rm Kod}(S/{i_3})=0$]
{$K^2=7,$ $g=3,$\\
${\rm Kod}(S/{i_1})=2,$ ${\rm Kod}(S/{i_2})=1,$ ${\rm Kod}(S/{i_3})=0$}\label{expl9}

This section contains the construction of a a bidouble cover
$V\rightarrow X,$ with $X$ rational, such that the minimal model
$S$ of $V$ is a surface of general type with $K^2=7,$ $p_g=q=1$ and $g=3$
having birational bicanonical map.

Let $i_j,$ $j=1,2,3,$ be the involutions associated to the bidouble cover.
The quotients $W_j:=V/{i_j},$ $j=1,2,3,$ satisfy:
\begin{description}
\item[$\cdot$] ${\rm Kod}(W_1)=2,$ $p_g(W_1)=1,$ $q(W_1)=0;$
\item[$\cdot$] ${\rm Kod}(W_2)=1,$ $p_g(W_2)=0,$ $q(W_2)=0;$
\item[$\cdot$] ${\rm Kod}(W_3)=0,$ $p_g(W_3)=0,$ $q(W_3)=1.$
\end{description}

This is an example for Theorems \ref{prop0}, b), \ref{prop1}, b) and \ref{prp2}, a).

We keep Notation \ref{notation}.
\\\\
{\bf$\cdot$ Construction of $S$}
\\
From Appendix \ref{CuM}, there exist a curve $C$ of type $7(3,(2,2)_T^4,3)$
({\it i.e.} $C$ is a plane curve of degree $7$ with triple points at $p_0,p_5$ and a tacnode
at $p_i$ tangent to the line $T_i$ through $p_0,p_i,$ $i=1,\ldots,4$)
and a curve $G$ of type $6(2,(2,2)_T^4,1),$ both through points $p_0,\ldots,p_5,$
such that $C+G$ is reduced.

Set:
$$
\begin{array}{l}
D_1:=\widetilde{T_1}+\widetilde{T_2}+\widetilde{T_3}-\sum_1^3 2E_i'+E_5-3E_0,\\
D_2:=\widetilde{T_4}+\widetilde{G}-\sum_1^3(2E_i+2E_i')-(2E_4+4E_4')-E_5-3E_0,\\
D_3:=\widetilde{C}-\sum_1^4(2E_i+2E_i')-3E_5-3E_0.
\end{array}
$$
We have
$$
\begin{array}{l}
L_1\equiv 7\widetilde T-\sum_1^3(2E_i+2E_i')-(2E_4+3E_4')-2E_5-3E_0,\\
L_2\equiv 5\widetilde T-\sum_1^3(E_i+2E_i')-(E_4+E_4')-E_5-3E_0,\\
L_3\equiv 5\widetilde T-\sum_1^4(E_i+2E_i')-3E_0
\end{array}
$$
and
$$
\begin{array}{l}
K_X+L_1\equiv 4\widetilde T-\sum_1^3(E_i+E_i')-(E_4+2E_4')-E_5-2E_0,\\
K_X+L_2\equiv 2\widetilde T-\sum_1^3 E_i'-2E_0,\\
K_X+L_3\equiv 2\widetilde T-\sum_1^4 E_i'+E_5-2E_0.
\end{array}
$$
Let $\psi:V\rightarrow X$ be the bidouble cover determined by $D_1,$ $D_2,$ $D_3$
and $S$ be the minimal model of $V.$ Then
$$\chi(\mathcal O_S)=4+\frac{1}{2}\sum_1^3 L_i(K_X+L_i)=4+0-1-2=1.$$
The Magma procedure $LinSys$ defined in \cite{Ri2} can be used to verify that
$$h^0(X,\mathcal O_X(K_X+L_1))=1$$ and
$$\begin{array}{l}
h^0(X,\mathcal O_X(2K_X+L_1+L_2))=1,\\
h^0(X,\mathcal O_X(2K_X+L_1+L_3))=1,\\
h^0(X,\mathcal O_X(2K_X+L_2+L_3))=1,
\end{array}
$$
therefore
$$p_g(S)=\sum_1^3h^0(X,\mathcal O_X(K_X+L_i))=1+0+0=1$$ and
the bicanonical map of $V$ is not composed with any of the involutions $i_1,$ $i_2,$ $i_3.$

The {\bf Albanese fibration} of $S$
is given by the pullback of the pencil of curves of type $6(2,(2,2)_T^4).$
\\\\
{\bf$\cdot$ Calculation of $K_S^2$}
\\
We have
$$N:=2K_X+\sum_1^3L_i\equiv N_1+N_2,$$
where $$N_1:=\widetilde{T_1}+\cdots+\widetilde{T_4}-4E_0-\sum_1^4 2E_i',\ \ \ N_2:=7\widetilde T-3E_0-\sum_1^4(2E_i+2E_i')-E_5.$$
As the support of the pullback of $N_1$ is a disjoint union of $16=-N_1^2$ $(-1)$-curves,
$|N_2|$ has no fixed component (use the procedure $LinSys$) and $N_1N_2=0,$
$$K_S^2=N^2+(-N_1^2)=N_2^2=7.$$
\\\\
{\bf$\cdot$ Verification that $\phi_2$ is birational}
\\
The system $|\psi^*(N)|$ is strictly contained in the bicanonical system of $V.$
The bicanonical map of $V$ is
not composed with any of the involutions
$i_1,$ $i_2,$ $i_3,$ hence it is birational if the map
$\tau$ given by $|N|=N_1+|N_2|$ is birational.
This is in fact the case, see Appendix \ref{CuM},
where Magma is used to show that the image of $\tau$ is of degree
$7=N_2^2.$
\subsection[$K^2=6,$ $g=3,$ ${\rm Kod}(S/{i_1})=2,$
${\rm Kod}(S/{i_2})=1,$ ${\rm Kod}(S/{i_3})=0$]
{$K^2=6,$ $g=3,$\\
${\rm Kod}(S/{i_1})=2,$ ${\rm Kod}(S/{i_2})=1,$ ${\rm Kod}(S/{i_3})=0$}\label{expl10}

One can obtain a construction analogous to the one in Section \ref{expl9}, but with $K_S^2=6$ instead:
replace the triple point of $C$ by a $(2,2)$-point, tangent to $G.$ Such a curve exists,
see Appendix \ref{CuM}.
With this change the branch locus in $W_3$ has a $4$-uple point instead of a $(3,3)$-point.
\subsection[$K^2=8,$ $g=3,$ ${\rm Kod}(S/{i_1})=1,$ $S/{i_2}$ ruled, ${\rm Kod}(S/{i_3})=1$]
{$K^2=8,$ $g=3,$\\
${\rm Kod}(S/{i_1})=1,$ $S/{i_2}$ ruled, ${\rm Kod}(S/{i_3})=1$}\label{expl7}

Here we give the construction of a surface of
general type $S,$ with $K^2=8,$ $p_g=q=1$ and $g=3,$
as a bidouble cover of a ruled surface $Z$ with $q(Z)=1.$

Let $i_j,$ $j=1,2,3,$ be the involutions associated to the bidouble cover.
The quotients $W_j:=S/{i_j},$ $j=1,2,3,$ satisfy:
\begin{description}
\item[$\cdot$] ${\rm Kod}(W_1)=1,$ $p_g(W_1)=0,$ $q(W_1)=1;$
\item[$\cdot$] $W_2$ is ruled, $q(W_2)=1;$
\item[$\cdot$] ${\rm Kod}(W_3)=p_g(W_3)=q(W_3)=1.$
\end{description}

This is an example for cases c) and d) of Theorem \ref{prop1}.

Notation \ref{notation} is used in this section.
\\\\
{\bf$\cdot$ Construction of $S$}
\\
Let $F_1,\ldots,F_4$ be disjoint fibres of the Hirzebruch surface $\mathbb F_0$
and $Z\rightarrow \mathbb F_0$ be the double cover with branch locus
$F_1+\cdots+F_4.$ Clearly $Z$ is a ruled surface with irregularity 1.
Denote by $\gamma$ the rational fibration of $Z.$

Let $G,G_1,\ldots,G_6$ be distinct smooth elliptic sections of $\gamma$ and
$\Gamma_1,\ldots,\Gamma_4$ be distinct fibres of $\gamma$ such
that $\Gamma_1+\Gamma_2\equiv 2\Gamma_3\equiv 2\Gamma_4.$

Set
$$
\begin{array}{l}
D_1:=\Gamma_1+\Gamma_2,\\
D_2:=G_1+\cdots+G_4,\\
D_3:=G_5+G_6
\end{array}
$$
and
$$
\begin{array}{l}
L_1:=3G+\Gamma_3-\Gamma_4,\\
L_2:=G+\Gamma_4,\\
L_3:=2G+\Gamma_3.
\end{array}
$$
We have
$$
\begin{array}{l}
K_Z+L_1\equiv G+\Gamma_3-\Gamma_4,\\
K_Z+L_2\equiv -G+\Gamma_4,\\
K_Z+L_3\equiv \Gamma_3
\end{array}
$$
and then
$$p_g(S)=\sum_1^3h^0(Z,\mathcal O_Z(K_Z+L_i))=0+0+1=1,$$
$$\chi(\mathcal O_S)=4\cdot 0+\frac{1}{2}\sum_1^3 L_i(K_Z+L_i)=0+0+0+1=1.$$

As
$$\begin{array}{l}
h^0(X,\mathcal O_X(2K_X+L_1+L_2))=0,\\
h^0(X,\mathcal O_X(2K_X+L_1+L_3))=2,\\
h^0(X,\mathcal O_X(2K_X+L_2+L_3))=1,
\end{array}
$$
the bicanonical map of $V$ is not composed with any of the involutions $i_1,$ $i_2,$ $i_3.$

The linear system
$$|N|:=|2K_Z+\sum_1^3L_i|\equiv |2G+2\Gamma_3|$$
has no base component, thus $$K_S^2=N^2=8.$$

\subsection[$K^2=4,$ $g=3,$ $S/{i_1}$ ruled, ${\rm Kod}(S/{i_2})=0,$ ${\rm Kod}(S/{i_3})=2$]
{$K^2=4,$ $g=3,$\\
$S/{i_1}$ ruled, ${\rm Kod}(S/{i_2})=0,$ ${\rm Kod}(S/{i_3})=2$}\label{expl11}

This section contains the construction of a bidouble cover
$V\rightarrow Z,$ with $Z$ ruled and $q(Z)=1,$ such that the minimal
model $S$ of $V$ is a surface of general type with
$K^2=4,$ $p_g=q=1,$ $g=3$ and that the bicanonical map
$\phi_2$ of $S$ is not composed with any of the
involutions $i_1,i_2,i_3$ associated to the bidouble cover.

The quotients $W_j:=V/{i_j},$ $j=1,2,3,$ satisfy:
\begin{description}
\item[$\cdot$] $W_1$ is ruled, $q(W_1)=1;$
\item[$\cdot$] ${\rm Kod}(W_2)=0,$ $p_g(W_2)=0,$ $q(W_2)=1;$
\item[$\cdot$] ${\rm Kod}(W_3)=2,$ $p_g(W_3)=1,$ $q(W_3)=1;$
the branch locus of the cover $V\rightarrow W_3$ is an union of four $(-2)$-curves.
\end{description}

This is an example for Theorems \ref{prop0}, b) and \ref{prp2}, b).

We use Notation \ref{notation}.
\\\\
{\bf$\cdot$ Construction of $S$}
\\
Let $Q_1$ be a general curve of type $5(1,(2,2)^3_T)$
(there is a pencil of such curves, see Appendix \ref{UsefulPencils})
and $Q_2$ be a general curve of type $3(1,(1,1)^3_T),$ both through points $p_0,\ldots,p_3.$

Let
$$Q_1':=\widetilde{Q_1}-\sum_1^3(2E_i+2E_i')-E_0\equiv 5\widetilde T-\sum_1^3(2E_i+2E_i')-E_0,$$
$$Q_2':=\widetilde{Q_2}-\sum_1^3(E_i+E_i')-E_0\equiv 3\widetilde T-\sum_1^3(E_i+E_i')-E_0.$$
Consider the double cover $\psi:Z\rightarrow X$ with branch locus
$$\widetilde{T_1}+\cdots+\widetilde{T_4}-\sum_1^3 2E_i'-4E_0,$$
where $T_4$ is a general line through $p_0.$

Let $$\Gamma:=\frac{1}{2}\psi^*(\widetilde{T_4}-E_0),\ \Gamma_i:=\frac{1}{2}\psi^*(\widetilde{T_i}-E_0),$$
$$C_0:=\psi^*(E_0),$$
$$e_i:=\frac{1}{2}\psi^*(E_i-E_i'),$$
$$e_i':=\psi^*(E_i')\ \ i=1,2,3,$$
and set
$$
\begin{array}{l}
D_1:=\psi^*(Q_1')\equiv 4C_0+10\Gamma-\sum_1^3(4e_i+4e_i'),\\
D_2:=\psi^*(Q_2')\equiv 2C_0+6\Gamma-\sum_1^3(2e_i+2e_i'),\\
D_3:=0,
\end{array}
$$
$$
\begin{array}{l}
L_1:=C_0+3\Gamma-\sum_1^3(e_i+e_i'),\\
L_2:=2C_0+5\Gamma-\sum_1^3(2e_i+2e_i'),\\
L_3:=3C_0+8\Gamma-\sum_1^3(3e_i+3e_i').
\end{array}
$$
Let $V$ be the bidouble cover of $Z$ determined by the curves $D_i$
and by the divisors $L_i$ and let $S$ be the minimal model of $V.$
Since $Z$ is a ruled surface with $q=1,$
$$K_Z\equiv -2C_0-2\Gamma+\sum_1^3\left[(e_i+e_i')+e_i\right].$$
Therefore,
$$
\begin{array}{l}
K_Z+L_1=-C_0+\Gamma+\sum_1^3 e_i,\\
K_Z+L_2=3\Gamma-\sum_1^3 e_i',\\
K_Z+L_3=C_0+6\Gamma-\sum_1^3(e_i+2e_i').
\end{array}
$$
One has
$$\chi(\mathcal O_S)=4\cdot 0+\frac{1}{2}\sum_1^3 L_i(K_Z+L_i)=0+0+1=1.$$
Let $$G=\widetilde{T_1}+\widetilde{T_2}+\widetilde{T_3}-\sum_1^3(E_i+E_i')+E_4-2E_0$$
and $$l=2\widetilde T-\sum_1^3E_i'-2E_0.$$
Since the branch locus of $\psi$ is linearly equivalent to $2l,$
$K_Z+L_3\equiv\psi^*(G)$ and $G-l\equiv\widetilde T-\sum_1^3E_i,$ then
$$h^0(Z,\mathcal O_Z(K_Z+L_3))=h^0(X,\mathcal O_X(G))+h^0(X,\mathcal O_X(G-l))=1+0=1.$$
This way,
$$p_g(S)=\sum_1^3h^0(Z,\mathcal O_Z(K_Z+L_i))=0+0+1=1.$$

As
$$\begin{array}{l}
h^0(Z,\mathcal O_Z(2K_Z+L_1+L_2))=0,\\
h^0(Z,\mathcal O_Z(2K_Z+L_1+L_3))>0,\\
h^0(Z,\mathcal O_Z(2K_Z+L_2+L_3))>0,
\end{array}
$$
the bicanonical map of $V$ is not composed with any of the involutions $i_1,$ $i_2,$ $i_3,$
corresponding to $L_1,L_2,L_3.$
\\\\
{\bf$\cdot$ Calculation of $K_S^2$}
\\
We have
$$N:=2K_Z+\sum_1^3L_i\equiv N_1+N_2,$$
where $$N_1:=2\Gamma_1+2\Gamma_2+2\Gamma_3-\sum_1^3 2e_i'$$
and
$$N_2:=2C_0+6\Gamma-\sum_1^3 (2e_i+2e_i')\equiv \psi^*\left(3\widetilde T-\sum_1^3(E_i+E_i')-E_0\right).$$
As the support of the pullback of $N_1$ is a disjoint union of $24=-N_1^2$ $(-1)$-curves,
$|N_2|$ has no fixed component and $N_1N_2=0,$ then
$$K_S^2=N^2+(-N_1^2)=N_2^2=4.$$

\appendix
\section[Appendix: Construction of plane curves]{Appendix: Construction of plane curves}\label{usefulcurves}

\subsection{Useful pencils}\label{UsefulPencils}

Here we show the existence of some pencils of plane curves that are useful on
some of the constructions of Section \ref{bdce}.
Recall Notation \ref{notation}.

\begin{lemma}\label{LemmaCbcCnc}
Let $C\subset\mathbb P^2$ be a smooth conic and $p_0\notin C,$ $p_1,\ldots,p_4\in C$
be distinct points. Consider the points $p_5,p_6\in C$ such that
the lines through $p_0,p_5$ and $p_0,p_6$ are tangent to $C.$

There exists a smooth curve $Q$ of type $3(1,(1,1)^4_T,1^2),$ through $p_0,\ldots,p_6.$
\end{lemma}
{\bf Proof:}
Let $C_x,$ $x\in\mathbb P^1,$ be a parametrization of the pencil of conics through $p_1,\ldots,p_4.$
Let $p_x^1,$ $p_x^2$ be the points of $C_x$ (not distinct if $C_x$ is singular) such that the
lines through $p_0,p_x^1$ and $p_0,p_x^2$ are tangent to $C_x.$
The correspondence $$x\leftrightarrow \{p_x^1,p_x^2\}$$ gives a plane algebraic curve $Q,$
parametrized by $x\in\mathbb P^1,$ and a double cover $Q\rightarrow\mathbb P^1.$ This cover is ramified over
four points, corresponding to the three degenerate conics which contain the points $p_1,\ldots,p_4$
plus the conic which contains $p_0.$ Therefore, by the Hurwitz formula, $Q$ is a cubic.

The conic through $p_0,\ldots,p_4$  is not tangent to the line $T_i$ (through $p_0,p_i$) at $p_0,$
thus also $Q$ is not tangent to $T_i$ at $p_0,$ $i=1,\ldots,4.$
Since each conic $C_x$ can be tangent to $T_i$ only at $p_i,$ $i=1,\ldots,4,$ then $Q$ intersects $T_i$
only at $p_0$ and $p_i,$ $i=1,\ldots,4.$ This means that $Q$ is tangent to $T_i$ at $p_i,$ $i=1,\ldots,4,$
and then $Q$ is smooth.\Qed
\begin{proposition}\label{PropUsefulPencils}
In the notation of Notation {\rm \ref{notation}},
there exist pencils, without base components, of
plane curves of type:
\begin{description}
  \item[\cite{Br} a)] $5(1,(2,2)^3_T);$
  \item[b)] $6(2,(2,2)^4_T);$
  \item[c)] $7(3,(2,2)^5_T);$
  \item[d)] $8(4,(2,2)^6_T).$
\end{description}
\end{proposition}
{\bf Proof:}
\begin{description}
  \item[a)] This is proved in \cite{Br}. Notice that we are imposing 19 conditions to a linear
  system of dimension 20.
  \item[b)] Let $\mathbb A(\mathbb C)$ be an affine plane and $a,b,c,d\in\mathbb C\backslash\{0\}$
  be numbers such that $a\ne c$ and $bc\ne\pm ad.$
  Consider the points of
  $\mathbb A:$ $$p_0:=(0,0),\ p_1:=(a,b),\ p_2:=(c,d),\ p_3:=(c,-d), p_4:=(a,-b)$$
  and let $T_i$ be the line through $p_0$ and $p_i,$ $i=1,\ldots,4.$
  Let $C_1$ be the conic through $p_1,\ldots,p_4$ tangent to $T_1,T_4$ and
  $C_2$ be the conic through $p_1,\ldots,p_4$ tangent to $T_2,T_3.$

  The curves $$2C_1+T_2+T_3\ \ \ {\rm and}\ \ \ 2C_2+T_1+T_4$$ generate a pencil
  whose general member is a curve of type $6(2,(2,2)^4_T).$
  \item[c)] Let $C\subset\mathbb P^2$ be a non-degenerate conic and $p_0\not\in C,$
  $p_1,\ldots,p_5\in C$
  be distinct points such that the lines $T_1,$ $T_5,$ defined by $p_0,p_1$ and $p_0,p_5,$ are
  tangent to $C.$ From Lemma \ref{LemmaCbcCnc}, there exists a curve $Q$
  of type $3(1,(1,1)^4_T,1),$ through $p_0,\ldots,p_5,$ respectively.

  The curves $$2C+T_2+T_3+T_4\ \ \ {\rm and}\ \ \ 2Q+T_5$$ generate a
  pencil whose general member is a
  curve of type $7(3,(2,2)^5_T).$
  \item[d)] This is analogous to the previous case, but now the pencil is generated by
  $$2C+T_2+\cdots+T_5\ \ \ {\rm and}\ \ \ 2Q+T_1+T_6.$$\Qed
\end{description}

\subsection{Constructions using Magma}\label{CuM}

In this appendix we construct some curves using the Computational Algebra System Magma (\cite{BCP}).
We use the Magma procedure $LinSys,$ defined in \cite{Ri2}.
This procedure calculates the linear system $L$ of plane curves of
degree $d,$ in an affine plane $\mathbb A,$ having singular points $p_i$ of
order $(m1_i,m2_i)$ with tangent direction given by the slope $td_i.$\\

Consider, in a affine plane $\mathbb A,$ the points
$$p_0:=(0,0),\ p_1:=(2,2),\ p_2:=(-2,2),\ p_3:=(3,1),\ p_4:=(-3,1).$$
From Section \ref{UsefulPencils}, there exists a pencil of curves of type
$6(2,(2,2)^4_T),$ with singularities at $p_0,\ldots,p_4,$ respectively.
Let $G$ be the element of this pencil which contains the point $p_5:=(0,5).$\\\\
\ \ $\cdot$ The curve $G$ is reduced and the tangent line to $G$ at $p_5$ is horizontal:
\begin{verbatim}
> A<x,y>:=AffineSpace(Rationals(),2);
> p:=[A![2,2],A![-2,2],A![3,1],A![-3,1],A![0,5],Origin(A)];
> d:=6;m1:=[2,2,2,2,1,2];m2:=[2,2,2,2,1];
> td:=[p[i][2]/p[i][1]:i in [1..4]] cat [0];
> LinSys(A,d,p,m1,m2,td,~L);
> #Sections(L);
1
> G:=Curve(A,Sections(L)[1]);
> IsReduced(G);
true
\end{verbatim}
\ \ $\cdot$ There exists a reduced curve $C$ of type $8(4,(2,2)_T^4,(3,3)),$
singular at $p_0,\ldots,p_5,$
such that the $(3,3)$-point is tangent to $G.$ Moreover, $G+C$ is a reduced element
of a pencil of curves of type $14(6,(4,4)_T^4,(4,4)):$
\begin{verbatim}
> d:=8;m1:=[2,2,2,2,3,4];m2:=[2,2,2,2,3];
> LinSys(A,d,p,m1,m2,td,~L);
> #Sections(L);
1
> C:=Curve(A,Sections(L)[1]);
> IsReduced(G join C);
true
> d:=14;m1:=[4,4,4,4,4,6];m2:=[4,4,4,4,4];
> LinSys(A,d,p,m1,m2,td,~L);
> #Sections(L);BaseComponent(L);
2
Scheme over Rational Field defined by
1
> #Sections(LinearSystem(L,G join C));
1
\end{verbatim}
Analogously, one can verify that:
\\\\
\ \ $\cdot$ there exists a reduced curve of type $7(3,(2,2)_T^4,(2,2)),$ singular at $p_0,\ldots,p_5,$
such that the $(2,2)$-point is tangent to $G.$
\\

Finally we will see that $p_5$ can be chosen such that\\\\
\ \ $\cdot$ there exist reduced curves $C_1$ of type $7(3,(2,2)_T^4,3)$ and $C_2$ of type\\
$6(2,(2,2)_T^4,1),$ both through $p_0,\ldots,p_5,$ such that $C_1+C_2$ is reduced and
the singularity of $C_1+C_2$ at $p_5$ is ordinary.

\begin{verbatim}
> A<x,y>:=AffineSpace(Rationals(),2);
> p:=[A![2,2],A![-2,2],A![3,1],A![-3,1],Origin(A)];
> d:=7;m1:=[2,2,2,2,3];m2:=[2,2,2,2];
> td:=[p[i][2]/p[i][1]:i in [1..#m2]];
> LinSys(A,d,p,m1,m2,td,~L);
> #Sections(L);BaseComponent(L);
6
Scheme over Rational Field defined by
1
\end{verbatim}
Now we impose a triple point to the elements of L.
This is done by asking for the vanishing of minors of a matrix of derivatives.
\begin{verbatim}
> R<x,y,n>:=PolynomialRing(Rationals(),3);
> h:=hom<PolynomialRing(L)->R|[x,y]>;
> H:=h(Sections(L));
> M:=[[H[i],D(H[i],1),D(H[i],2),D2(H[i],1,1),D2(H[i],1,2),\
> D2(H[i],2,2)]:i in [1..#H]];
> Mt:=Matrix(M);min:=Minors(Mt,#H);
> A:=AffineSpace(R);
> S:=Scheme(A,min cat [x-3,1+n*(y-x)*(y+x)*(3*y-x)*(3*y+x)]);
> //The condition 1+n*(..)=0 implies that
> //the solution is not in p.
> Dimension(S);
0
> PointsOverSplittingField(S);
\end{verbatim}
We choose one of the solutions and show that it works:
\begin{verbatim}
> R<r1>:=PolynomialRing(Rationals());
> K<r1>:=NumberField(r1^2 - 1761803/139426560*r1 + \
> 1387488001/33730073395200);
> A<x,y>:=AffineSpace(K,2);
> y1:=-33462374400/102856069*r1 + 419793163/102856069;
> p:=[A![2,2],A![-2,2],A![3,1],A![-3,1],A![3,y1],Origin(A)];
> d:=7;m1:=[2,2,2,2,3,3];m2:=[2,2,2,2];
> td:=[p[i][2]/p[i][1]:i in [1..#m2]];
> LinSys(A,d,p,m1,m2,td,~L);#Sections(L);
1
> C1:=Curve(A,Sections(L)[1]);
> IsOrdinarySingularity(C1,p[5]);
true
> d:=6;m1:=[2,2,2,2,1,2];m2:=[2,2,2,2];
> LinSys(A,d,p,m1,m2,td,~L);#Sections(L);
1
> C2:=Curve(A,Sections(L)[1]);
> IsReduced(C1 join C2);
true
> IsSingular(C2,p[5]);
false
> IsOrdinarySingularity(C1 join C2,p[5]);
true
\end{verbatim}
The verification that the singularities are no worst than stated is left to the reader
(use the Magma functions {\em ProjectiveClosure}, {\em SingularPoints}, {\em HasSingularPointsOverExtension}
and {\em ResolutionGraph}).\\

The calculations for Section \ref{expl9} (verification that $\phi_2$ is birational) are as follows:
\begin{verbatim}
> d:=7;m1:=[2,2,2,2,1,3];m2:=[2,2,2,2];
> LinSys(A,d,p,m1,m2,td,~L);
> #Sections(L);BaseComponent(L);
5
Scheme over K defined by
1
> P4:=ProjectiveSpace(K,4);
> tau:=map<A->P4|Sections(L)>;
> Degree(tau(Scheme(A,Sections(L)[3])));
7
\end{verbatim}
thus an hyperplane section of the image of $\tau$ is of degree 7.

\bibliography{ReferencesRito3}

\bigskip
\bigskip

\noindent Carlos Rito
\\ Departamento de Matem\' atica
\\ Universidade de Tr\' as-os-Montes e Alto Douro
\\ 5000-911 Vila Real
\\ Portugal
\\\\
\noindent {\it e-mail:} crito@utad.pt

\end{document}